\documentclass[11pt,a4paper]{article}%

\usepackage{amsthm,amssymb}
\usepackage{enumitem}

\usepackage{authblk}
\theoremstyle{definition}

\usepackage[T1]{fontenc}
\usepackage{graphicx}

\usepackage{lscape}
\usepackage[utf8]{inputenc}
\usepackage{lmodern} %
\usepackage[english]{babel}
\usepackage{tabularx}
\usepackage{amsmath,amscd,latexsym}
\usepackage{bm}
\usepackage{hyperref}
\usepackage{subcaption}
\usepackage{algorithm}
\usepackage{algpseudocode}
\algrenewcommand\algorithmicensure{\textbf{Input:}}
\usepackage{color}
\usepackage{graphicx}
\graphicspath{{figures/}{tikzs/}{svg-inkscape/}}
\usepackage{overpic}
\usepackage{subcaption}
\usepackage{placeins} %
\usepackage{cleveref}
\usepackage{siunitx}
\usepackage{varwidth}

\usepackage{tikz}
\usetikzlibrary{shapes.geometric,calc,positioning}
\usepackage{tikz-cd}
\usepackage{tikz-3dplot}
\usepackage{pgfplots}				
\usepackage{varwidth}

\newlength\figurewidth
\newlength\figureheight
\newcommand{\includetikz}[1]{\includegraphics{tikzs/#1_compiled}}

\let\originalleft\left
\let\originalright\right
\renewcommand{\left}{\mathopen{}\mathclose\bgroup\originalleft}
\renewcommand{\right}{\aftergroup\egroup\originalright}

\newcommand{\tunderbrace}[2]{\underbrace{#1}_{\textstyle#2}}

\usepackage[normalem]{ulem}
\usepackage{graphicx}
\usepackage{wrapfig}

\newcommand{\displacement}{w}
\newcommand{\const}{\textrm{const.}}

\begin{document}
\title{On an optimization framework for damage localization in structures\thanks{This work is part of dtec.bw -- Digitalization and Technology Research Center of the
		Bundeswehr within project “Digitalisierung von Infrastrukturbauwerken zur Bau\-werks\-über\-wachung: Structural Health Monitoring”. dtec.bw is funded by the European Union~-- NextGenerationEU.}}
\author[1]{Owais Saleem\thanks{\texttt{saleemo@hsu-hh.de}}}
\author[1]{Tim Suchan\thanks{\texttt{suchan@hsu-hh.de}}}
\author[1]{Natalie Rauter\thanks{\texttt{natalie.rauter@hsu-hh.de}}}
\author[1]{Kathrin Welker\thanks{\texttt{welker@hsu-hh.de}}}
\affil[]{Helmut Schmidt University/University of the Federal Armed Forces Hamburg, Department of Mechanical and Civil Engineering, Holstenhofweg 85, 22043 Hamburg, Germany}
\date{}

\providecommand{\keywords}[1]
{
  \small	
  \textbf{\textit{Key words. }} #1
}
\providecommand{\ams}[1]
{
  \small	
  \textbf{\textit{AMS subject classifications. }} #1
}
\maketitle              %
	\begin{abstract}
	Efficient structural damage localization remains a challenge in structural health monitoring (SHM), particularly when the problem is coupled with uncertainty of conditions and complexity of structures. Traditional methods simply based on experimental data processing are often not sufficiently reliable, while complex models often struggle with computational inefficiency given the tremendous amount of model parameters.
	This paper focuses on closing the gap between data-driven SHM and physics-based model updating by offering a solution for real-world infrastructure. We first concentrate on fusing multi-source damage-sensitive features (DSF) based on experimental modal data into spatially mapped belief masses to pre-screen candidate damage locations. The resulting candidate damage locations are integrated into an inverse Finite Element method (FEM) model calibration
	process. We propose an optimization framework to identify the most probable damage scenario with single and multi-damage cases. We present the corresponding numerical results in this paper, which open the door to extend the application of the framework to a complex real bridge structure.
	\end{abstract}

\section{Introduction} \label{sec:Introduction}
	Structural health monitoring aims at detecting damage, deterioration or abnormal behavior. It has rapidly evolved into a cornerstone of Industry~4.0 for critical infrastructure, employing advanced sensor networks, real\-time big data analytics, and digital twins to enable truly predictive maintenance.  A classical SHM lifecycle---often referred to as the Rytter scheme \cite{rytter1993}---comprises four stages: detection, localization, quantification, and prognosis.  Our hybrid FEM–experimental modal analysis (EMA) framework focuses squarely on the second and third stages, i.e., on accurately pinpointing damage sites and assessing their severity.
		
	Industry~4.0 has shifted SHM from passive inspections into active intelligence systems, in which edge‑computing devices preprocess streaming measurements and cloud platforms run high‑fidelity simulations in parallel.  In this context, damage localization and quantification becomes the linchpin: accurate spatial mapping of defects not only triggers focused maintenance actions, but also refines the digital twin for continuous model updating. %
	
	Damage localization in SHM has been approached from both model-based and data-driven perspectives, and there are modal-static hybrid and uncertainty-quantification schemes:
	\begin{itemize}
	\item \emph{Model-based approaches}: FEM-based model updating formulates damage localisation as an inverse problem in which selected model parameters (e.g. element stiffness reductions, local flexibility or boundary properties) are adjusted to minimise discrepancies between measured and simulated dynamic observables (natural frequencies, mode shapes, frequency response functions or FRFs) \cite{Friswell1995}.
	Sensitivity-based modal updating and least-squares formulations provide efficient gradient information for Gauss–Newton or Levenberg–Marquardt solvers in modal and FRF objectives \cite{Mottershead2011}.
	FRF-based and hybrid FRF/modal updating schemes have been shown to improve localisation sensitivity to local stiffness changes in beam studies \cite{Zhan2021}.
	Practical two-step and substructural updating strategies have also been demonstrated for laboratory and field-scale validation of damage-localisation procedures \cite{Malekghaini2023}.
	\item \emph{Data-driven methods}:  Autoregressive and statistical pattern-recognition algorithms detect damage by monitoring shifts in time-series model coefficients. In \cite{Lu2005}, an Auto-Regressive with exogenous inputs (ARX)-based prediction model whose residual-error statistics serve as damage indicators is proposed. Vector auto-regressive moving-average (VARMA) and stochastic subspace methods for damage detection and localization under unknown excitation are reviewed in \cite{Bodeux2003}.
	Moreover, empirical-mode decomposition is combined with VARMA to define a damage index in operational environments (cf.~\cite{Dong2010}). Evidence‐based fusion via Dempster–Shafer theory has been shown to be effective for structural damage localization in beam‐type systems. First, \cite{Fei2009} demonstrated---through a simulation study---how damage sensitivity theory can combine modal residuals to pinpoint damage in an Euler‐Bernoulli beam, achieving high localization accuracy under varying noise levels.  Then, \cite{Bao2011} extends this framework by formulating basic belief assignments from mode‐shape and frequency changes and applies Dempster’s rule to detect damage in lab‐scale structures, which shows improved robustness compared to traditional probabilistic methods. 
	In addition, \cite{Li2008} further enhanced damage sensitivity-based identification by integrating Shannon‐entropy-weighted evidence fusion, which adaptively discounts conflicting or low‐information features, yielding more reliable damage indices in both simulation and experimental validations.
	\item \emph{Modal–static hybrids and uncertainty-aware schemes:} Modal-curvature approaches compare experimentally measured and analytically computed curvature fields to highlight local stiffness changes and produce damage indicators. In \cite{Farrar2007}, the fusion of changes in low-order modal parameters with higher-order or complementary information to improve localization and diagnosis is discussed. A framework for assessing the robustness of SHM technologies that emphasizes principle feature selection, classifier choice and uncertainty-informed decision making is presented in \cite{Stull2012}, thus motivating the fusion of multiple feature types (e.g., static, dynamic or statistical) under noisy conditions. Finally, \cite{Kefal2016} introduces a four-node quadrilateral inverse-shell element with drilling degrees of freedom for iFEM/shape sensing, demonstrating improved full-field reconstruction, which is useful for sensing and damage assessment in plates and composite panels.
	\end{itemize}
	
	Both the power and the limitations of a standalone inverse FEM-based damage quantification framework is discussed and demonstrated and a novel framework is introduced. This new framework concentrates on fusing multi-source damage-sensitive features based on experimental modal data into spatially mapped belief masses to pre-screen candidate damage locations. The resulting candidate damage locations are integrated into an inverse Finite Element model calibration process. 
	
	Even though our approach provides a general framework for structures, we consider beam-like structures that extend only in one dimension for simplicity reasons, with the goal of detecting local damage. The previously-mentioned simulations will be using an Euler-Bernoulli beam model, which is described by the partial differential equation with the spatial coordinate $x$ and time $t$, i.e.,
	\begin{align}
		\rho(x) A(x) \frac{\partial^2 \displacement(x,t)}{{\partial t}^2} + \frac{\partial^2}{{\partial x}^2} \left( E(x) I(x) \frac{\partial^2 \displacement(x,t)}{{\partial x}^2} \right) = q(x,t).
		\label{eqn:EulerBernoulliBeamEquation}
	\end{align}
	Here, $\rho(x)$ denotes the density of the considered material, $A(x)$ is the cross-sectional area, $\displacement(x)$ is the beam deflection, $E(x)$ describes the Young's modulus, $I(x)$ is the second moment of area of the beam's cross section, and $q(x,t)$ describes external loading. The partial differential equation also requires boundary and initial conditions, which depend on the specific setting that is simulated. Damage is then described by damage parameters $\bm{\theta}$---common choices are Young's modulus, density, or a combination of both. It may be noted that out of the two candidates for damage modelling (Young's modulus and density), it is well established that a crack significantly influences the former while mass distribution remains virtually unchanged, cf., e.g.,~\cite{Cawley1979}.

	In this paper, we propose a hybrid framework which addresses Rytter stages~2 (localization) and~3 (quantification), fusing full-field reconstructions with modal indicators under a Dempster-Shafer filter to yield both precise damage maps and severity estimates.
	
	The paper is organized as follows:
	The fundamentals of Dempster-Shafer theory is explained in Section~\ref{sec:DS Localization}. Here, two sources of damage are considered: modal strain energy change ratio (MSECR) and a flexibility-based index. The two sources are fused and a more condensed plausibility information is utilized for localization.
	Section~\ref{sec:iFEM_Localization} defines our objective for the optimization process and discusses the optimization strategy along with the chosen algorithm. The damage quantification approach for purely FEM-based optimization is discussed, along with its limitations, and two new strategies are formulated and implemented: the hierarchical approach and the coupled DS-FEM approach.
	In Section~\ref{sec:results}, we present numerical results for each implemented strategy: evi\-dence-theo\-retic localization, pure FEM-based localization and quantification, and the fully coupled Dempster-Shafer and FEM framework.
	Finally, Section~\ref{sec:conclusion} summarizes the findings, discusses practical implications for SHM, and outlines directions for extending the hybrid framework.

\section{Damage localization based on evidence theory}
\label{sec:DS Localization}

Within this work, the Dempster-Shafer (DS) evidence theory \cite{Dempster67,Shafer76} is used in the context of damage detection. DS evidence theory extends classical Bayesian inference by permitting belief mass to be assigned not only to single hypotheses but to sets of hypotheses, thereby explicitly representing ignorance when evidence is inconclusive. The theory is briefly introduced in the following.

	\subsection{Fundamental definitions}

Let the frame of discernment
\begin{align*}
\Theta = \{H_1,\dots,H_n\}
\end{align*}
be the set of all possible mutually exclusive outcomes or hypotheses for a given situation, where the power set \(2^\Theta\) contains all subsets \(A\subseteq\Theta\), including the empty set. In this context a basic probability assignment (BPA) or mass function
\begin{align*}
m\colon2^\Theta\to[0,1] 
\end{align*}
with the properties
\begin{align*}
\sum_{A\subseteq\Theta}m(A)=1
\quad \text{and} \quad
m(\emptyset)=0,
\end{align*}
assigns the exact belief mass assigned to each subset \(A\) and no other subset \cite{Shafer76}. One defines the belief and plausibility from the mass \(m\) of any \(A\subseteq\Theta\) as
\begin{align*}
	\mathrm{Bel}(A) &= \sum_{B\subseteq A} m(B), \\
	\mathrm{Pl}(A) &= \sum_{B\cap A\neq\emptyset} m(B).
\end{align*}
Moreover, if \(\overline{A}=\Theta\setminus A\) denotes the complement of \(A\) in the frame \(\Theta\), then the plausibility can be written as
\begin{align*}
	\mathrm{Pl}(A)
	&= \sum_{B\cap A\neq\emptyset} m(B)
	= 1 - \sum_{B\subseteq\overline{A}} m(B)
	= 1 - \mathrm{Bel}\left(\overline{A} \right).
\end{align*}
Hence, we obtain \(\mathrm{Bel}(A)\le \mathrm{Pl}(A)\), and the interval \([\,\mathrm{Bel}(A),\ \mathrm{Pl}(A)\,]\subseteq[0,1]\) quantifies the total and the maximally plausible support for \(A\) given the available evidence.
		
	\paragraph{Dempster’s rule of combination.}  
		Given two independent BPA \(m_1,m_2\) on the same \(\Theta\), their orthogonal sum \(m=m_1\oplus m_2\) is defined by (cf.~\cite{Dempster67}):
		\begin{align*}
			m(C)  =  \frac{1}{1 - K}
			\sum_{A\cap B = C}
			m_1(A)\,m_2(B),
			\quad
			K = \sum_{A\cap B=\emptyset} 
			m_1(A)\,m_2(B),
		\end{align*}
To compute \(m(C)\), first sum the products \(m_1(A)m_2(B)\) over all subsets \(A\) and \(B\) whose intersection equals \(C\). Then, normalize this sum by dividing by \(1-K\), where \(K\) quantifies the total conflict between \(m_1\) and \(m_2\) as the mass assigned to disjoint subsets.
		
\subsection{Dynamic ignorance and feature-specific BPA mapping} \label{sec:dynamic_ignorance}

To map damage-sensitive features $\bm{D} = \left( D_1,\dots,D_n \right)^\top$ into BPA that explicitly quantify uncertainty, we introduce a physics-grounded dynamic ignorance model combining four evidence cues. This approach addresses key limitations of static ignorance models in structural health monitoring \cite{Basir2007,SentzFerson02}.

\subsubsection{Normalization and concentration metrics}
		
For each element $i$ (where each element is a section of the beam between uniformly distributed measurement points) with damage index $D_i \geq 0$, compute:
\begin{align*}
D^{\text{norm}}_i &= \frac{D_i}{\sum_{k=1}^n D_k+ \epsilon}  & \text{(probability mass)} \\
D^{\text{rel}}_i &= \frac{D_i}{\max_k D_k + \epsilon} & \text{(relative severity)} 
\end{align*}
Here, $\epsilon = 10^{-10}$ prevents division by zero. 
			
Using the probability mass, we can then quantify the damage concentration via normalized Shannon entropy \cite{Shannon48}:
\begin{align}
H(\bm{D}) = -\frac{1}{\ln n} \sum_{k=1}^n D^{\text{norm}}_k \ln (D^{\text{norm}}_k + \epsilon), \quad
{conc} = 1 - H(\bm{D})
\label{eq:entropy}
\end{align}
Here, ${conc} \approx 1$ indicates localized damage (low uncertainty), while ${conc} \approx 0$ suggests distributed damage (high uncertainty).

\subsubsection{Ranking and ignorance components} \label{ssec:ignorance_components}

The dynamic ignorance model decomposes uncertainty into four physically interpretable components by following established frameworks that distinguish between different types of ignorance and uncertainty in decision-making \cite{smithson2022ambiguities}.
			
\begin{figure}[tbp]
\centering
\setlength\figurewidth{0.6\textwidth}%
\setlength\figureheight{6cm}%
\includetikz{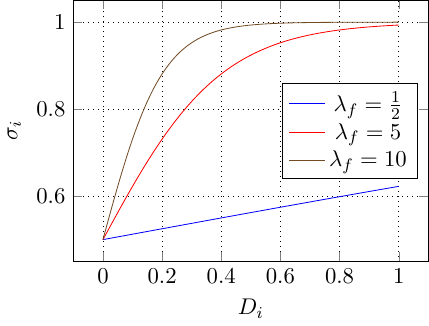}%
\caption{Feature‐specific logistic confidence curves.}
\label{fig:confidence_curve}
\end{figure}

\begin{enumerate}
	\item The \textbf{distribution uniformity }
		\begin{align*}
		d_i = 1 - {conc}
		\end{align*}
	measures global damage dispersion uncertainty using information entropy for each element $i$. High values ($d_i \approx 1$) occur when damage is uniformly distributed (\mbox{$H(\bm{D}) \approx 1$}), indicating low localization confidence. Low values ($d_i \approx 0$) indicate concentrated damage with high confidence \cite{Shannon48}. This aligns with the concept of complexity-based uncertainty where system properties emerge from multiple interacting elements \cite{taleb2012antifragile}.
				
	\item For each element $i$, the \textbf{relative weakness}
		\begin{align*}
		r'_i = 1 - D^{\text{rel}}_i
		\end{align*}
	quantifies element-specific damage severity uncertainty using relative damage indices, approaches 1 for elements with $D_i \ll \max(\bm{D})$ (low relative damage), and approaches 0 for maximally damaged elements. This component reflects probability-based uncertainty where outcomes are indeterminate but estimable through relative comparisons \cite{rovelli2020probability}.
				
	\item The \textbf{rank-based uncertainty}
		\begin{align*}
		k_i = \frac{r_i}{n-1}
		\end{align*}
	with \(r_i\in\{0,\dots,n-1\}\), where \(0\) denotes the highest damage, represents positional uncertainty in the damage hierarchy for each element $i$. Thus, top-ranked elements have $k_i=0$ and bottom-ranked elements have $k_i=1$. This addresses ambiguity-based uncertainty where information reliability varies across rankings \cite{han2019relationship}.
				\item For each element $i$, the \textbf{confidence-based uncertainty}
				\begin{align*}
					c_i = 1 - \sigma_i, \quad \text{where} \quad \sigma_i = \frac{1}{1 + \exp(-\lambda_f D_i)}
				\end{align*}
				characterizes measurement reliability using the feature-specific sensitivity $\lambda_f$ by incorporating the psychological concept of \emph{tolerance for uncertainty} \cite{han2022perceptions}. This logistic mapping transforms damage indices into confidence scores, reflecting how different features contribute to ignorance \cite{ferson1996different}. Each feature type $f$ (frequency, strain, etc.) has an associated sensitivity parameter $\lambda_f > 0$.
				Higher values of $\lambda_f$ produce steeper confidence transitions, cf. \Cref{fig:confidence_curve}.  
				Assigning a larger $\lambda_f$ to a particular feature~$f$ means expressing greater trust in that feature's damage‐index values—i.e., even a modest non‑zero $D_i$ from that feature is strong evidence of real damage.
\end{enumerate}

\subsubsection{Weighted ignorance synthesis} \label{ssec:ignorance_synthesis}
The ignorance factor $\alpha_i$ combines the four ignorance and uncertainty components through physics-informed weighting, i.e.,
\begin{align}
\alpha_i = w_d d_i + w_r r'_i + w_k k_i + w_c c_i
\label{eq:total_ignorance}
\end{align}
with weights $\bm{w} = \left( w_d, w_r, w_k, w_c \right)^\top \in (0,1)^4$ related to the distribution uniformity $d_i$, the relative weakness $r'_i$, the rank-based uncertainty $k_i$, and the feature confidence $c_i$, respectively. The weights are constrained by $\sum w = 1$.
More details about the weights and their physical significance are provided in Table \ref{tab:weight_interpretation}. 
			
			\begin{table}[tbp]
				\caption{Weight interpretation and physical significance in ignorance synthesis}
				\centering
				\begin{tabularx}{\textwidth}{|l|X|}
					\hline
					Weight & Effect when emphasized \\
					\hline
					$w_d$ (distribution) &
					\begin{itemize}[leftmargin=*,nosep]
						\item Suppresses uncertainty when damage is concentrated
					\end{itemize} \\
					\hline
					$w_r$ (relative) &
					\begin{itemize}[leftmargin=*,nosep]
						\item Amplifies uncertainty for marginally damaged elements
					\end{itemize} \\
					\hline
					$w_k$ (rank) &
					\begin{itemize}[leftmargin=*,nosep]
						\item Prioritizes damage hierarchy over magnitude
					\end{itemize} \\
					\hline
					$w_c$ (confidence) &
					\begin{itemize}[leftmargin=*,nosep]
						\item Boosts confidence for reliable features
					\end{itemize} \\
					\hline
				\end{tabularx}
				\label{tab:weight_interpretation}
			\end{table}
			
			The ignorance factor is constrained to $\alpha_i \in [\alpha_{\min}, 1]$ where $\alpha_{\min} = 0.1$, which ensures a minimum 10\% uncertainty even for strong evidence. We set $\alpha_{\min}=0.1$ (a 10\% ignorance floor) as a conservative lower bound on uncertainty to prevent the fusion process becoming over-confident when evidence is limited or when model/sensor errors and unmodelled effects are present \cite{yang2013evidential}.

		\subsubsection{Basic probability assignment (BPA)}	
		The basic probability assignment  for each element $i$ is defined by allocating mass to the singleton hypothesis $\{i\}$ and to the frame of discernment $\Theta$ (representing ignorance). The belief factor $\beta_i = 1 - \alpha_i$ scales the normalized damage index, i.e.,
		\begin{align*}
			m_i(\{i\}) &= \beta_i \cdot D^{\text{norm}}_i \text{ and } \\
			m_i(\Theta) &= \alpha_i,
		\end{align*}
		where
		\begin{itemize}
			\item $D^{\text{norm}}_i\colon$ describes the normalized damage index for element $i$,
			\item $\beta_i\colon$ indicates the belief factor, quantifying the specific belief committed to element $i$, and
			\item $\alpha_i\colon$ represents the ignorance factor, quantifying the uncommitted belief mass.
		\end{itemize}
	
\subsection{Modal strain‐energy change ratio (MSECR)}
		\label{ssec: ds_msecr}
		The modal strain‐energy change ratio (MSECR) quantifies the relative change in strain‐energy that is stored in each Finite Element between undamaged and damaged states \cite{Fei2009}. In the Finite Element discretization the classical elastic strain energy associated with a (static) displacement vector \(\bm{u}\) is \(\frac{1}{2}\bm{u}^\top \bm{K}\bm{u}\). Here, we evaluate the corresponding modal contribution by inserting the modal displacements (eigenvectors). Let $\{\,\bm{\phi}_j^{und}\,\}$ and $\{\,\bm{\phi}_j^{dam}\,\}$ be the mass‐normalized eigenvectors (more details in \Cref{ssec:objective_theory}) of the undamaged and damaged beam, respectively, for modes $j=1,\dots,m$. Denote  the elemental stiffness matrix for element $i$ by $\bm{K}_i$. All elemental stiffness matrices \(\mathbf{K}_i\) are assembled in the global degree-of-freedom basis and therefore share the same dimension as the global eigenvector \(\boldsymbol{\phi}\).
		
		We first compute the total modal strain‐energy in element $i$:
		\begin{align*}
		{SE}_i^{und}
		= \sum_{j=1}^m \left( \bm{\phi}_j^{und} \right)^\top \bm{K}_i \, \bm{\phi}_j^{und},
		\quad
		SE_i^{dam}
		= \sum_{j=1}^m \left( \bm{\phi}_j^{dam} \right)^\top \bm{K}_i\,\bm{\phi}_j^{dam}.
		\end{align*}
		The MSECR index is then defined as the absolute relative change:
		\begin{align*}
			D_i^{SE}
			= \left| SE_i^{und} - SE_i^{dam} \right|
		\end{align*}
		Elements experiencing a local stiffness reduction (e.g., due to damage) exhibit a pronounced drop in modal strain‐energy, making $D_i^{SE}$ a sensitive indicator of damage localization. 
		
		As is clear from the formulation, the MSECR relies on both the experimental data as well as element stiffness matrices, which are retrieved from a FEM solver.

	\subsection{Flexibility‐based damage index}
		\label{ssec: ds_flexibility}
		
		Following Grande and Imbimbo’s multi‐stage damage detection framework \cite{Grande2007}, we construct an element‐wise flexibility change index that directly measures how each element’s resistance to bending decreases due to damage. Let the global flexibility (inverse of stiffness) matrices for the undamaged and damaged structures be obtained by modal superposition:
		\begin{align*}
		\bm{F}^{und} = \sum_{j=1}^m \frac{1}{\left( \omega_j^{und} \right)^2} \bm{\phi}_j^{und} \left( \bm{\phi}_j^{und} \right)^\top,
		\quad
		\bm{F}^{dam} = \sum_{j=1}^m \frac{1}{\left( \omega_j^{dam} \right)^2} \bm{\phi}_j^{dam} \left( \bm{\phi}_j^{dam} \right)^\top
		\end{align*}
		For each element $e$, we extract its local curvature‐strain compliance via the projection matrices $\bm{L}_e$ (disassembly) and $\bm{S}_e$ (curvature extraction):
		\begin{align*}
		\bm{F}_e^{\mathrm(\cdot)}
		= \bm{S}_e \left( \bm{L}_e \bm{F}^{(\cdot)} {\bm{L}_e}^\top \right) {\bm{S}_e}^\top,
		\quad
		(\cdot)\in\{und,dam\}
		\end{align*}
		The element index then is the Frobenius norm of the compliance change, i.e.,
		\begin{align*}
			D_e^{F}
			= \left\| \bm{F}_e^{dam} - \bm{F}_e^{und}\right\|_F,
		\end{align*}
		which quantifies the increase in curvature-per-unit-moment (i.e., curvature compliance) at element $e$, corresponding to a local reduction in flexural rigidity $EI$ (as $\kappa/M=1/(EI)$). This flexibility‐based measure has been shown to be highly sensitive to local stiffness reductions, making it an effective complement to modal‐energy‐based indices \cite{Grande2007}.
		
		We have now briefly introduced Dempster-Shafer evidence theory. In the next section, we continue with damage localization using FEM.
	
\section{FEM based model updating for damage localization}
	\label{sec:iFEM_Localization}
	
	In this section, we describe an approach of damage detection and localization based on optimization, in combination with a numerical (Finite Element) model. In \Cref{ssec:objective_theory} we first introduce the objective functional that is to be minimized. \Cref{ssec:OptimizationStrategy} introduces a gradient-based optimization approach, together with the required derivatives. In \Cref{ssec:hierarchical} we describe an additional hierarchical process that uses the optimization approach and iteratively locates damage in a more and more detailed manner.
	
	\subsection{Objective function}
	\label{ssec:objective_theory}
	
	We formulate the damage‐localization problem as an optimization problem in which the total objective function $J$ is the weighted sum of three error terms. The objective function reads	
	\begin{align}
		J(\bm{\theta}) = \alpha_f E_{f}(\bm{\theta})+\alpha_g E_{g}(\bm{\theta})+\alpha_c E_{c}(\bm{\theta})\,,
		\label{eqn:ObjectiveFunctional}
	\end{align}
	where $\bm{\theta}$ denotes the vector of damage parameters and $\alpha_{(\cdot)}$ are the corresponding weights for each part of the objective function. The three contributions in \eqref{eqn:ObjectiveFunctional} serve complementary roles in the objective:
	\begin{itemize}
		\item $E_f\colon$ Frequency shifts detect global stiffness changes, making this term sensitive to overall structural integrity. 
		\item $E_g\colon$ The governing residual enforces model-experiment consistency, thereby acting as a physics-based objective. 
		\item $E_c\colon$ Curvature captures local bending anomalies, providing good spatial resolution for damage localization.
	\end{itemize}
	
		For the following descriptions, we consider a classical free vibration Euler-Bernoulli beam as in \eqref{eqn:EulerBernoulliBeamEquation}, i.e., $q(x,t)=0$. We use the common assumption that the system response (i.e., the time-dependent deflection of the beam $w(x,t)$) can be decomposed into a time- and a space-dependent part, i.e.,
		\begin{align*}
			w(x,t) = X(x) \cdot T(t), \text{ where } T(t) = \mathrm{e}^{-\mathrm{i} \omega t}.
		\end{align*}
		This yields the differential equation
		\begin{align}
			\begin{aligned}
			0 
			&= \rho(x) A(x) X(x) \frac{\partial^2 T(t)}{{\partial t}^2} + \frac{\partial^2}{{\partial x}^2} \left( E(x) I(x) \frac{\partial^2 X(x)}{{\partial x}^2} \right) T(t) \\
			&= \rho(x) A(x) X(x) (-\mathrm{i} \omega)^2 \mathrm{e}^{-\mathrm{i} \omega t} + \frac{\partial^2}{{\partial x}^2} \left( E(x) I(x) \frac{\partial^2 X(x)}{{\partial x}^2} \right) \mathrm{e}^{-i \omega t} \\
			&= -\rho(x) A(x) X(x) \, \omega^2  + \frac{\partial^2}{{\partial x}^2} \left( E(x) I(x) \frac{\partial^2 X(x)}{{\partial x}^2} \right).
			\end{aligned}
			\label{eqn:EulerBernoulliBeamEquationODE}
		\end{align}
		If we have constant material parameters $E(x)=E=\const$, $I(x)=I=\const$, $\rho(x)=\rho=\const$ and $A(x)=A=\const$, then we immediately obtain
		\begin{align*}
			0 = -\rho A \omega^2 X(x) + E I \frac{\partial^4 X(x)}{{\partial x}^4}.
		\end{align*}
		Solving this ordinary differential equation for $X(x)$ yields infinitely many solutions $X_j(x)$, $j=1,2,\ldots$, which are given by
		\begin{align*}
			X_j(x) = A \cosh(\beta_j x) + B \sinh(\beta_j x) + C \cos(\beta_j x) + D \sin(\beta_j x),
		\end{align*}
		where the constants $A$, $B$, $C$ and $D$ depend on the boundary conditions of the differential equation~\eqref{eqn:EulerBernoulliBeamEquation}. These boundary conditions prescribe the value of $X_j(x)$ or its derivatives at some points $x\in \left[ -\frac{L}{2}, \frac{L}{2} \right]$, where $L$ is the length of the beam. The natural frequencies $\omega_j$ and $\beta_j$ are related to each other by
		\begin{align*}
			\beta_j = \sqrt[4]{\frac{\rho A \omega_j^2}{E I}}.
		\end{align*}
		Common descriptions of this problem also include the eigenvalue $\lambda_j$, which simply is given by $\lambda_j = \omega_j^2$. The beam vibration frequency can also easily be calculated as $f_j = \frac{\omega_j}{2\pi}$. The corresponding function $X_j(x)$ is called the mode and the shape of $X_j(x)$, denoted as $\phi_{j}(x)$, is called the mode shape or eigenvector.

		However, non-constant material parameters, e.g., due to damage, require alternative methods to obtain frequencies and mode shapes. The Finite Element method is commonly used to obtain (a approximation of) these frequencies and mode shapes. This is the approach that we follow in this manuscript. Thus, our mode shapes and frequencies are obtained numerically and in a discrete sense. We use a uniform computational mesh with element size~$h$. Only the first $m$ modes are considered, as higher-frequency modes tend to be more difficult to measure in physical experiments due to higher structural damping, smaller amplitudes, and thus a stronger influence of noise on the measurements. As already described in \Cref{sec:Introduction}, damage in the form of a crack is best modeled by a local change in the Young's modulus. Thus, we identify the damage parameters~$\bm{\theta}$ with the Young's modulus~$E(x)$.

	\paragraph{Frequency‐shift error.}
	We measure relative shifts of the first $m$ eigenfrequencies between experiment and model. The experimental data, and therefore also the eigenfrequencies from the experiment, are independent of the choice of (numerical) damage parameters $\bm{\theta}$. The model eigenfrequency has a dependence on $\bm{\theta}$. Thus, we have
	\begin{align*}
		E_{f}(\bm{\theta}) = \sum_{j=1}^{m} 
		\left(\frac{\omega_{j}^{exp}-\omega_{j}^{mod}(\bm{\theta})}{\omega_{j}^{exp}}\right)^{2}.
	\end{align*}
	
	\paragraph{Governing‐equation residual.}
	If the numerical model matches perfectly with the experiment, then all experimental mode shapes $\bm{\phi}_{j}^{exp}(x)$ fulfill \eqref{eqn:EulerBernoulliBeamEquationODE}. Using the Finite Element method yields a discretization, i.e., a discretized mode shape~$\bm{\phi}_j^{mod}(x,\bm{\theta})$ (a vector), a global mass matrix~$\bm{M}^{mod}(\bm{\theta})$, and a global stiffness matrix $\bm{K}^{mod}(\bm{\theta})$. Then, f%
	or each mode $j$, the ideal relation
	\begin{align*}
		\left( \omega_{j}^{exp} \right)^{2}\,\left( \bm{\phi}_{j}^{exp}\right)^\top \bm{M}^{mod}(\bm{\theta}) \,\bm{\phi}_{j}^{exp} = \left( \bm{\phi}_{j}^{exp}\right)^\top \bm{K}^{mod}(\bm{\theta}) \,\bm{\phi}_{j}^{exp}
	\end{align*}
	should hold. By reordering, we define the normalized residual error
	\begin{align*}
		E_{g}(\bm{\theta}) =
		\sum_{j=1}^{m}
		\left(1 - 
		\frac{\left( \bm{\phi}_{j}^{exp}\right)^\top \bm{K}^{mod}(\bm{\theta}) \,\bm{\phi}_{j}^{exp}}
		{\left( \omega_{j}^{exp} \right)^{2}\,\left( \bm{\phi}_{j}^{exp}\right)^\top \bm{M}^{mod}(\bm{\theta}) \,\bm{\phi}_{j}^{exp}}
		\right)^{2}.
	\end{align*}
	
	\paragraph{Curvature‐error term.} The curvature~$\kappa$ for each mode $j$ at a point $x$ along the beam can be defined as the change of the tangent vector, with a deflection according to the mode shape $\bm{\phi}_{j}$. Thus, the curvature of mode $j$ can be expressed as
	\begin{align*}
		\kappa_j(x,\bm{\theta}) = \frac{\frac{\partial^2}{{\partial x}^2} \left( \phi_{j}(x,\bm{\theta}) \right)}{\left(1 + \left(\frac{\partial}{\partial x} \left( \phi_{j}(x,\bm{\theta}) \right)^2 \right)^{\frac{3}{2}} \right)}.
	\end{align*}
	As commonly encountered in beam theory, for small curvatures we obtain 
	\begin{align*}
		\left| \frac{\partial}{\partial x} \left( \phi_{j}(x,\bm{\theta}) \right) \right| \ll 1,
	\end{align*}
	and thus,
	\begin{align*}
		\kappa_j(x,\bm{\theta}) \approx \frac{\partial^2}{{\partial x}^2} \left( \phi_{j}(x,\bm{\theta}) \right).
	\end{align*}
	When considering the experimental curvature, no dependence on $\bm{\theta}$ is present, i.e., $\kappa_{j}^{exp}(x,\bm{\theta})=\kappa_{j}^{exp}(x)$.
	By comparing model and experimental curvatures for each mode $j$ at $N$ points $x_1, \ldots, x_N$, we obtain
	\begin{align*}
		E_{c}(\bm{\theta}) = \sum_{j=1}^{m}
		\sum_{i=1}^{N}
		\left(
		\frac{\kappa_{j}^{exp}(x_i)-\kappa_{j}^{mod}(x_i,\bm{\theta})}
		{\max\left( \left| \kappa_{j}^{exp}(x_i) \right|, \epsilon \right)}
		\right)^{2},
	\end{align*}
	where $0 < \epsilon \ll 1$ is a small regularization constant to avoid division by zero. The points $x_1, \ldots, x_N$, at which the curvature $\kappa_j$ is evaluated, are chosen as the (interior) nodes of the model.
	
	\subsection{Optimization strategy}
	\label{ssec:OptimizationStrategy}
	
	The damage parameters $\bm{\theta} = \left(\theta_{1}, \dots, \theta_{k} \right)^\top$ are estimated by minimizing $J(\bm{\theta})$ using the limited‐memory Broyden–Fletcher–Goldfarb–Shanno 
	(L--BFGS) algorithm with analytically computed gradients.  Let $\lambda_{j}^{mod}$ and $\bm{\phi}_{j}^{mod}$ be the $j$th eigenvalue and discretized eigenvector of the current model, as described in~\Cref{ssec:objective_theory}. The derivatives of $\lambda_{j}^{mod}$ and $\bm{\phi}_{j}^{mod}$ with respect to $\theta_k$, $k=1,\ldots,K$ are given in \cite{FoxKapoor} and read
	\begin{align*}
		&\frac{\partial \lambda_{j}^{mod}(\bm{\theta})}{\partial \theta_{k}}
		= {\bm{\phi}_{j}^{mod} (x,\bm{\theta})}^\top
		\left(
		\frac{\partial \bm{K}^{mod}(\bm{\theta})}{\partial \theta_{k}}
		-\lambda_{j}^{mod}(\bm{\theta}) \,\frac{\partial \bm{M}^{mod}(\bm{\theta})}{\partial \theta_{k}}
		\right)
		\bm{\phi}_{j}^{mod} (x,\bm{\theta}), \\
		&\begin{aligned}
		\frac{\partial \bm{\phi}_{j}^{mod} (x,\bm{\theta})}{\partial \theta_{k}}
		= \sum_{\substack{h = 1 \\ h \neq j}}^{H}
		&\frac{1}{\lambda_{j}^{mod} (\bm{\theta}) - \lambda_{h}^{mod} (\bm{\theta})} 
		\left( {\bm{\phi}_{h}^{mod} (x,\bm{\theta})}^\top \left(
		\frac{\partial \bm{K}^{mod} (\bm{\theta})}{\partial \theta_{k}}
	 	\right. \right. \\
		&\ \left. \left.-\lambda_{j}^{mod} (\bm{\theta})\,\frac{\partial \bm{M}^{mod} (\bm{\theta})}{\partial \theta_{k}} 
		\right)
		\bm{\phi}_{j}^{mod} (x,\bm{\theta})\right)
		\,\bm{\phi}_{h}^{mod} (x,\bm{\theta}),
		\end{aligned}
	\end{align*}
	where $\bm{M}^{mod} (\bm{\theta})$ and $\bm{K}^{mod} (\bm{\theta})$ are the global model mass and stiffness matrices, and \mbox{$H \leq m$} is the total number of retained modes (cf. \cite[Formulation~2]{FoxKapoor}). As previously mentioned, we only use the Young's modulus as an indication of damage and thus have $\frac{\partial \bm{M}^{mod} (\bm{\theta})}{\partial \theta_{k}}=\bm{0}$. The derivative of the local (element) stiffness matrix with respect to the Young's modulus of element~$i$ is easily calculated and is given by $\frac{1}{\theta_k} \cdot \bm{K}_{i}(\bm{\theta})$ if $i=k$, while the derivative with respect to $\theta_{k}$ with $k \neq i$ is zero. Reassembling the global stiffness matrix from these evaluations then yields $\frac{\partial \bm{K}^{mod} (\bm{\theta})}{\partial \theta_{i}}$. These expressions are used to form the gradient $\nabla J(\bm{\theta})$, which the L-BFGS optimizer then employs to update $\bm{\theta}$ iteratively until convergence.
	
	In order to perform optimization, we require derivatives of our objective with respect to our design variables, which are provided in the following. A more detailed calculation of these derivatives can be found in \Cref{app:CalculationDerivatives}.
	
	By using $\lambda_j^{mod}(\bm{\theta}) = \left( \omega_j^{mod} (\bm{\theta}) \right)^2$ the derivative of frequency-shift error reads
	\begin{align*}
		\begin{aligned}
		\frac{\partial E_f(\bm{\theta})}{\partial \theta_{k}}
		&= - \sum_{j=1}^{m} \frac{\omega_{j}^{exp}-\omega_{j}^{mod}(\bm{\theta})}{\left( \omega_{j}^{exp} \right)^2 \omega_{j}^{mod}(\bm{\theta}) } \cdot \frac{\partial \lambda_{j}^{mod}(\bm{\theta})}{\partial \theta_{k}}.
		\end{aligned}
	\end{align*}
	
	As the model mass matrix is independent of $\bm{\theta}$ in our case, the derivative of the Governing-equation residual is given by
	\begin{align*}
			&\frac{\partial E_g(\bm{\theta})}{\partial \theta_{k}} \\
			&= 2 \sum_{j=1}^{m} \left(1 - 
			\frac{\left( \bm{\phi}_{j}^{exp}\right)^\top \bm{K}^{mod}(\bm{\theta}) \,\bm{\phi}_{j}^{exp}}
			{\left( \omega_{j}^{exp} \right)^{2}\,\left( \bm{\phi}_{j}^{exp}\right)^\top \bm{M}^{mod} \,\bm{\phi}_{j}^{exp}}
			\right) \cdot \left( - 
			\frac{\left( \bm{\phi}_{j}^{exp}\right)^\top \frac{\partial \bm{K}^{mod}(\bm{\theta})}{\partial \theta_{k}} \,\bm{\phi}_{j}^{exp}}
			{\left( \omega_{j}^{exp} \right)^{2}\,\left( \bm{\phi}_{j}^{exp}\right)^\top \bm{M}^{mod} \,\bm{\phi}_{j}^{exp}} \right).\nonumber
	\end{align*}
	
	The derivative of the curvature-error term requires a finite difference approximation of the second derivative of the mode shape, i.e.,
	\begin{align*}
		\frac{\partial^2}{{\partial x}^2} \left( \phi_{j}^{mod}(x,\bm{\theta}) \right) \approx \bm{C} \begin{pmatrix}
			\phi_{j}^{mod}(x_1,\bm{\theta}) \\
			\phi_{j}^{mod}(x_2,\bm{\theta}) \\
			\vdots \\
			\phi_{j}^{mod}(x_N,\bm{\theta}) \\
		\end{pmatrix} = \bm{C} \bm{\phi}_{j}^{mod}(x,\bm{\theta}),
	\end{align*}
	where the matrix $\bm{C}$ is given by
	\begin{align*}
		\bm{C} = \frac{1}{h^2}
		\begin{pmatrix}
			2 & -5 & 4 & -1 &  \vphantom{\ddots} & 0 & 0 & 0 & 0 \\
			1 & -2 & 1 & 0 & \ldots\vphantom{\ddots} & 0 & 0 & 0 & 0 \\
			0 & 1 & -2 & 1 & \ddots & 0 & 0 & 0 & 0 \\
			0 & 0 & 1 & -2 & \ddots & 0 & 0 & 0 & 0 \\
			\vdots &\vdots &\ddots & \ddots & \ddots & \ddots & \ddots & \vdots & \vdots \\
			0 & 0 & 0 & 0 &  \ddots & -2 & 1 & 0 & 0 \\
			0 & 0 & 0 & 0 &  \ddots & 1 & -2 & 1 & 0 \\
			0 & 0 & 0 & 0 &  \ldots\vphantom{\ddots} & 0 & 1 & -2 & 1 \\
			0 & 0 & 0 & 0 & \vphantom{\ddots}  & -1 & 4 & -5 & 2
		\end{pmatrix}.
	\end{align*}

	An evaluation at $x=x_i$ for $i=1,\ldots,N$ yields the $i$-th entry of $\bm{C} \bm{\phi}_{j}^{mod}(x,\bm{\theta})$, and we obtain
	\begin{align*}
		\frac{\partial E_c(\bm{\theta})}{\partial \theta_{k}}  = -2 \sum_{j=1}^{m} \sum_{i=1}^{N} \frac{\kappa_{j}^{exp}(x_i)-\kappa_{j}^{mod}(x_i,\bm{\theta})}
		{\max\left( \left| \kappa_{j}^{exp}(x_i) \right|, \epsilon \right)^2} \cdot \left( \left. \bm{C} \frac{\partial \bm{\phi}_{j}^{mod}(x,\bm{\theta})}{\partial \theta_{k}} \right|_{x=x_i} \right).
	\end{align*}

	\paragraph{Challenges with inverse problems.}  
	Solving inverse problems via optimization is inherently ill‑posed: multiple parameter sets may yield similar model outputs (non‑uniqueness), small measurement errors can cause large fluctuations in the estimated parameters (instability) and yield erratic gradients as well.
	
	There is no direct solution to this problem as it is not the limitation of the algorithm, but the very ill‑posedness of the problem that stagnates the progress towards convergence to the true damage profile. However, there are ways to circumvent this problem by regularization or other means. In the context of FEM‑based damage localization, a popular technique is the Tikhonov regularization \cite{tikhonov1963} which has also been used in this study.
	
	\textbf{Tikhonov regularization:}
	To dampen oscillatory solutions and improve the stability of the inverse problem, a quadratic penalty on the damage vector is introduced, leading to the regularized objective function
	\begin{align*}
		J_{\text{reg}}(\bm{\theta}) = J(\bm{\theta}) + \gamma \left\| \bm{\theta} \right\|_{2}^{2},
	\end{align*}
	where $\gamma > 0$ is the regularization parameter. The choice of $\gamma$ is critical and can be guided by the L-curve criterion~\cite{Hansen1992}, which seeks a balance between the solution norm $\|\bm{\theta}\|_2$ and the residual norm $J(\bm{\theta})$.
	
	While Tikhonov regularization improves well-posedness, selecting an appropriate $\gamma$ and achieving a physically plausible solution can still be challenging. This motivated the exploration of the alternative solution strategies presented in the following sections.

	\subsection{Hierarchical approach}
	\label{ssec:hierarchical}
	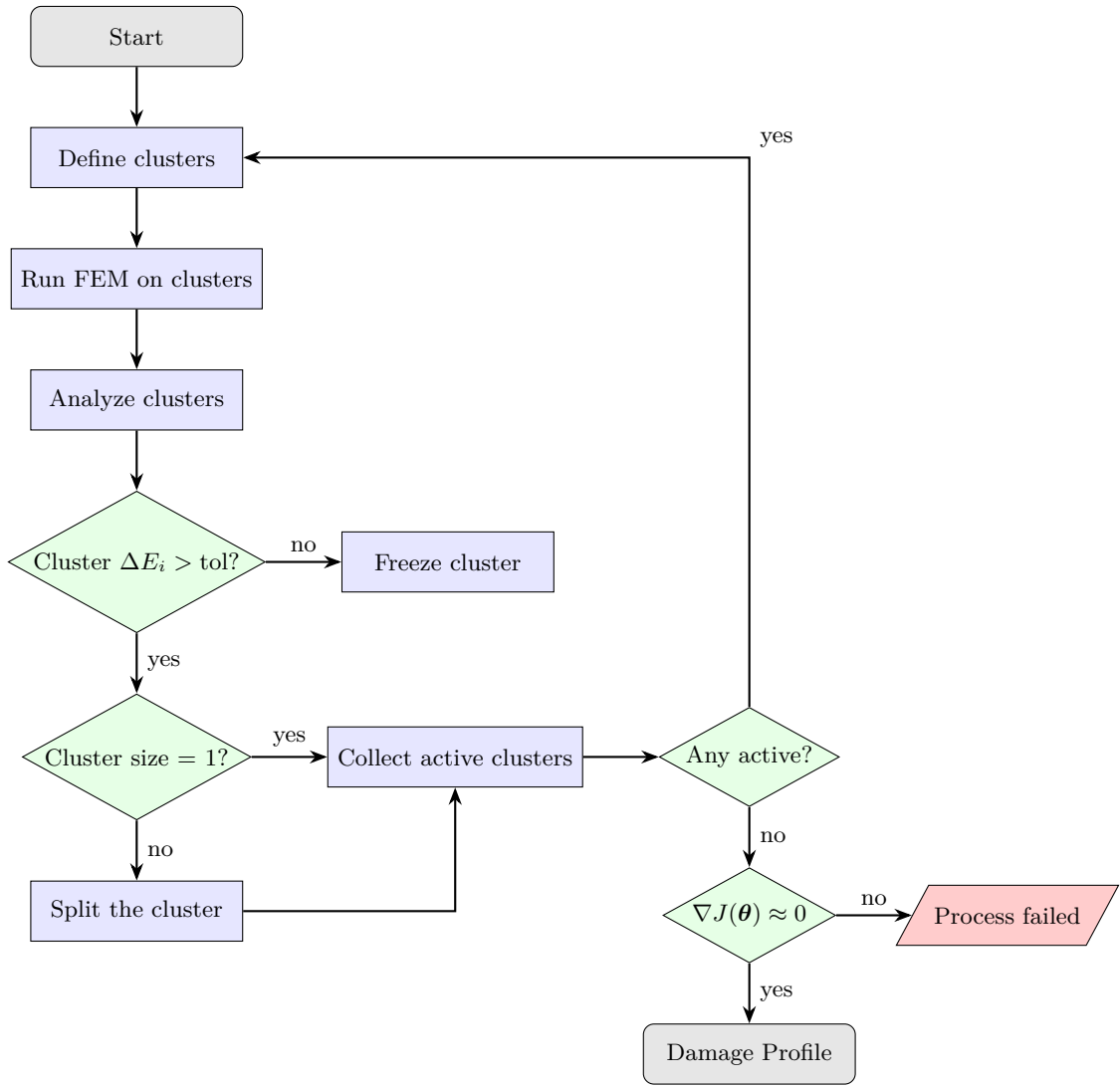
\begin{figure}[tbp]
		\centering
		
		\begin{tikzpicture}[
			node distance=8mm and 8mm,
			every node/.style={font=\footnotesize, align=center},
			startstop/.style={
				rectangle, rounded corners, draw, fill=gray!20,
				minimum width=2.8cm, minimum height=8mm
			},
			process/.style={
				rectangle, draw, fill=blue!10,
				minimum width=2.8cm, minimum height=8mm
			},
			decision/.style={
				diamond, draw, fill=green!10,
				aspect=1.8, inner sep=1pt, minimum width=1.8cm
			},
			io/.style={
				shape=trapezium,
				trapezium stretches=true,
				trapezium left angle=40,     %
				trapezium right angle=180-40,    %
				draw,
				fill=red!20,
				minimum width=2.8cm,
				minimum height=8mm,
				inner sep=2pt,
				align=center
			},
			flow/.style={-Stealth, thick},
			]
			
			\node (start)     [startstop]                  {Start};
			\node (define)    [process, below=of start]   {Define clusters};
			\node (select)    [process, below=of define]  {Run FEM on clusters};
			\node (analyze)   [process, below=of select]  {Analyze clusters};
			\node (processG)  [decision, below=of analyze] {Cluster $\Delta E_i > \mathrm{tol}$?};
			\node (checksize) [decision, below=of processG] {Cluster size = 1?};
			\node (split)     [process, below=of checksize]{Split the cluster};
			\node (freeze)    [process, right=of processG,xshift=2mm]   {Freeze cluster};
			\node (collect)   [process, right=of checksize,xshift=2mm]  {Collect active clusters};
			\node (checkact)  [decision, right=of collect,xshift=2mm]   {Any active?};
			\node (checkObj)  [decision, below=of checkact]   {$\nabla J(\bm{\theta})\approx 0 $};
			\node (message)   [io, right=of checkObj, xshift=2mm]     	{Process failed};
			\node (end)       [startstop, below=of checkObj]            {Damage Profile};
			\draw [flow] (start)    -- (define);
			\draw [flow] (define)   -- (select);
			\draw [flow] (select)   -- (analyze);
			\draw [flow] (analyze)  -- (processG);
			\draw [flow] (collect)  -- (checkact);
			\draw [flow] (processG.east) |- node[pos=0.75, above] {no} (freeze.west);
			\draw [flow] (checksize.south) -| node[pos=0.75, right] {no} (split.north);
			\draw [flow] (split.east) -| (collect.south);
			\draw [flow] (checkact.south) -| node[pos=0.75, right] {no} (checkObj.north);
			\draw [flow] (checkObj.east) |- node[pos=0.75, above] {no} (message.west);
			\draw[flow]	(processG.south) -|	node[pos=0.75,right] {yes}(checksize.north);
			\draw [flow] (checksize.east) |- node[pos=0.75, above] {yes} (collect.west);
			\draw [flow] (checkact.north) |- (define.east) node[midway, above right] {yes};
			\draw [flow] (checkObj.south) -| node[pos=0.75, right] {yes} (end.north);
			
		\end{tikzpicture}
		\caption{Hierarchical coarse‐to‐fine optimization framework}
		\label{fig:hierarchical_flowchart}
	\end{figure}

	To address the aforementioned challenges, one of the strategies implemented is a coarse-to-fine refinement strategy. We initiate the optimization with a minimal, coarse representation of the damage field. Upon identifying an optimal coarse damage pattern, we selectively refine regions flagged by the optimizer while freezing less sensitive clusters to their original (healthy) states. This process is sketched in \Cref{fig:hierarchical_flowchart}

	This progressive refinement constrains the solution space and mitigates optimizer stagnation, as each stage operates within a reduced subspace. Crucially, early iterations capture global damage trends, enabling more robust and efficient convergence during subsequent local refinements.
	
	The decision of how and when the optimizer refines or freezes clusters is governed by dynamic stage tolerances. Starting with the largest clusters, the process proceeds until stagnation occurs. The clusters are then evaluated based on their deviation from healthy values. The tolerance is defined as a specified percentage (e.g., 90\%) of the maximum observed variation. Clusters exceeding this tolerance are activated (refined) and values less than that are deemed insignificant changes with no physical meaning and hence are reset to their healthy values.

\section{Numerical results}
\label{sec:results}
	\begin{figure}[tbp]
		\centering
	\setlength\figurewidth{0.8\textwidth}%
	\newlength\svgwidth%
	\graphicspath{{figures/svg-inkscape/}}
	\setlength\svgwidth{\figurewidth}
\begingroup%
  \makeatletter%
  \providecommand\color[2][]{%
    \errmessage{(Inkscape) Color is used for the text in Inkscape, but the package 'color.sty' is not loaded}%
    \renewcommand\color[2][]{}%
  }%
  \providecommand\transparent[1]{%
    \errmessage{(Inkscape) Transparency is used (non-zero) for the text in Inkscape, but the package 'transparent.sty' is not loaded}%
    \renewcommand\transparent[1]{}%
  }%
  \providecommand\rotatebox[2]{#2}%
  \newcommand*\fsize{\dimexpr\f@size pt\relax}%
  \newcommand*\lineheight[1]{\fontsize{\fsize}{#1\fsize}\selectfont}%
  \ifx\svgwidth\undefined%
    \setlength{\unitlength}{270.61421792bp}%
    \ifx\svgscale\undefined%
      \relax%
    \else%
      \setlength{\unitlength}{\unitlength * \real{\svgscale}}%
    \fi%
  \else%
    \setlength{\unitlength}{\svgwidth}%
  \fi%
  \global\let\svgwidth\undefined%
  \global\let\svgscale\undefined%
  \makeatother%
  \begin{picture}(1,0.15064446)%
    \lineheight{1}%
    \setlength\tabcolsep{0pt}%
    \put(0,0){\includegraphics[width=\unitlength,page=1]{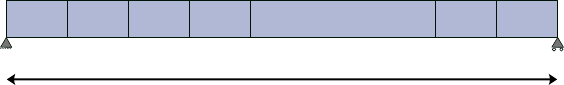}}%
    \put(0.41419533,0.02675253){\color[rgb]{0.14117647,0.12156863,0.10980392}\makebox(0,0)[lt]{\lineheight{1.25}\smash{\begin{tabular}[t]{l}1000 mm\end{tabular}}}}%
    \put(0.91316086,0.10474473){\color[rgb]{0,0,0.34117647}\makebox(0,0)[lt]{\lineheight{1.25}\smash{\begin{tabular}[t]{l}20\end{tabular}}}}%
    \put(0.80371997,0.10473647){\color[rgb]{0,0,0.34117647}\makebox(0,0)[lt]{\lineheight{1.25}\smash{\begin{tabular}[t]{l}19\\\end{tabular}}}}%
    \put(0.0757147,0.10470358){\color[rgb]{0,0,0.34117647}\makebox(0,0)[rt]{\lineheight{1.25}\smash{\begin{tabular}[t]{r}1\end{tabular}}}}%
    \put(0.1627733,0.10448984){\color[rgb]{0,0,0.34117647}\makebox(0,0)[lt]{\lineheight{1.25}\smash{\begin{tabular}[t]{l}2\\\end{tabular}}}}%
    \put(0.27136516,0.10474473){\color[rgb]{0,0,0.34117647}\makebox(0,0)[lt]{\lineheight{1.25}\smash{\begin{tabular}[t]{l}3\\\\\end{tabular}}}}%
    \put(0.38012352,0.10474473){\color[rgb]{0,0,0.34117647}\makebox(0,0)[lt]{\lineheight{1.25}\smash{\begin{tabular}[t]{l}4\\\end{tabular}}}}%
    \put(0,0){\includegraphics[width=\unitlength,page=2]{beam_svg-tex.pdf}}%
  \end{picture}%
\endgroup%
		\caption{Schematic of the one-dimensional beam.}
		\label{fig:beamSketch}
	\end{figure}
	All simulations are performed on a beam (cf. \Cref{fig:beamSketch}) with rectangular cross-section, which is 1000\,mm long, 20\,mm wide and 3.25\,mm thick. One-dimensional Euler-Bernoulli beam theory is adopted. Different boundary conditions are used such as simply-supported and cantilever. By default, results correspond to the simply-supported beam unless otherwise specified.
	
	The beam was discretized into 20 elements (21 nodes), which is a relatively fine discretization compared to, e.g., \cite{Fei2009}.
	
	\subsection{DS-based damage localization}

	Initially, the following is evaluated as the potential sources of damage information:
		\begin{enumerate}
			\item Damage indices based on eigenfrequencies
			
			Damage alters stiffness and thus natural frequencies.  For each element $i$, we compute the weighted sum of absolute frequency shifts $\Delta\omega_j=\left|\omega^{dam}_j-\omega^{und}_j \right|$, using normalized stiffness matrix sensitivities to form weights., i.e.,
			\begin{align*}
			D^f_i  =  \sum_{j=1}^{m} 
			\frac{\left|(\bm{\phi}^{dam}_j)^T \bm{K}_i\,\bm{\phi}^{dam}_j\right|}
				{\sum_{i=1}^{n_e}\left|(\bm{\phi}^{dam}_j)^T \bm{K}_i\,\bm{\phi}^{dam}_j\right|} 
			\left|\omega^{dam}_j-\omega^{und}_j\right|.
			\end{align*}

		\item Damage indices based on modal curvature
		
		For each mode \(j\) we calculate the curvature residuals \(\kappa^{und}_{j}-\kappa^{dam}_{j}\).
		Averaging squared residuals over \(m\) modes gives a nodal root-mean-square (RMS) \(R_n\), where \(n\) indexes nodal locations; this is projected to element \(i\) by averaging its end-node values:
		\[
		D^\kappa_i = \frac{1}{2}\bigl(R_{n_i} + R_{n_{i+1}}\bigr),\qquad
		R_n = \sqrt{\frac{1}{m}\sum_{j=1}^{m}
			\bigl(\kappa^{und}_{j}(x_n) - \kappa^{dam}_{j}(x_n)\bigr)^{2}}.
		\]
		We associate an element curvature to element \(i\) (which spans from node \(x_i\) to node \(x_{i+1}\)) by averaging its end-node curvatures:
		\[
		\kappa_{j}([x_i, x_{i+1}]) \approx \tfrac{1}{2}\bigl(\kappa_{j}(x_i) + \kappa_{j}(x_{i+1})\bigr).
		\]

			\item Damage indices based on strain energy
			
			Already described in detail in \Cref{ssec: ds_msecr}.
			\item Damage indices based on flexibility
			
			Already described in detail in \Cref{ssec: ds_flexibility}.
		\end{enumerate}
		
		\begin{figure}[tbp]
			\centering
			\begin{subfigure}[b]{0.48\textwidth}
	\setlength\figurewidth{\textwidth}%
	\newlength\svgwidth%
	\graphicspath{{figures/DS/svg-inkscape/}}
	\setlength\svgwidth{\figurewidth}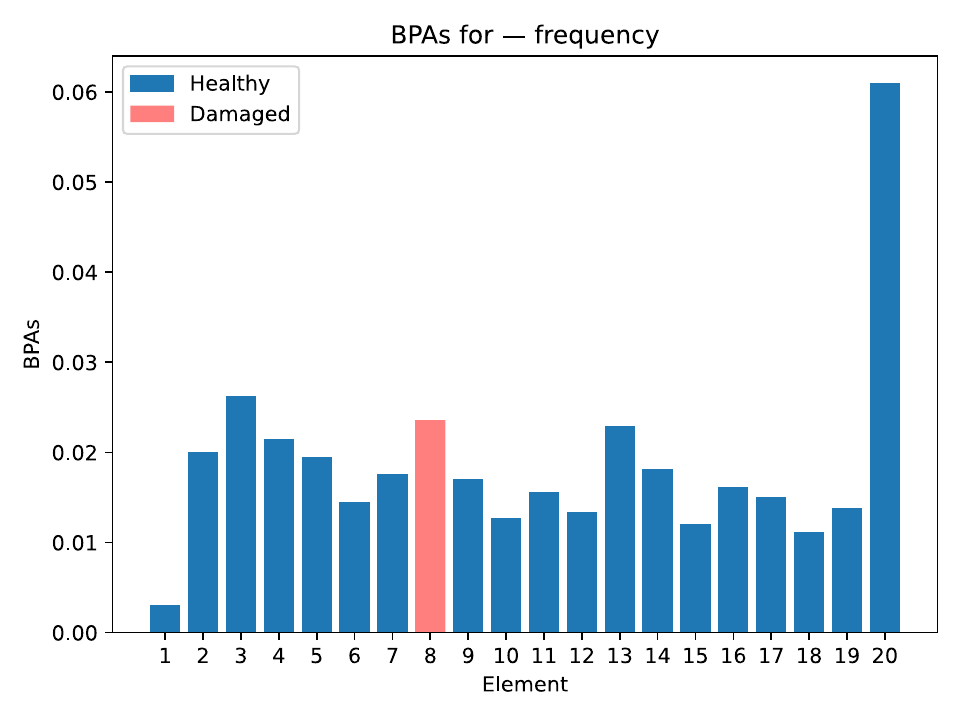%

				\label{fig:sub1}
			\end{subfigure}
			\hfill
			\begin{subfigure}[b]{0.48\textwidth}
	\setlength\figurewidth{\textwidth}%
	\newlength\svgwidth%
	\graphicspath{{figures/DS/svg-inkscape/}}
	\setlength\svgwidth{\figurewidth}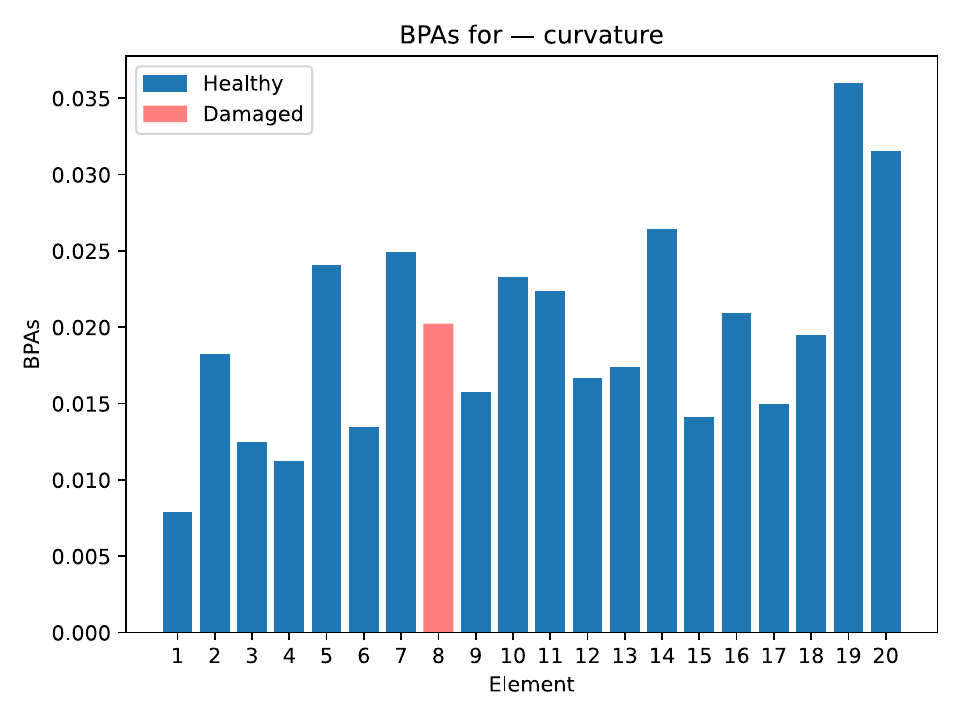%

				\label{fig:sub2}
			\end{subfigure}
			\vspace{1em}
			\begin{subfigure}[b]{0.48\textwidth}
	\setlength\figurewidth{\textwidth}%
	\newlength\svgwidth%
	\graphicspath{{figures/DS/svg-inkscape/}}
	\setlength\svgwidth{\figurewidth}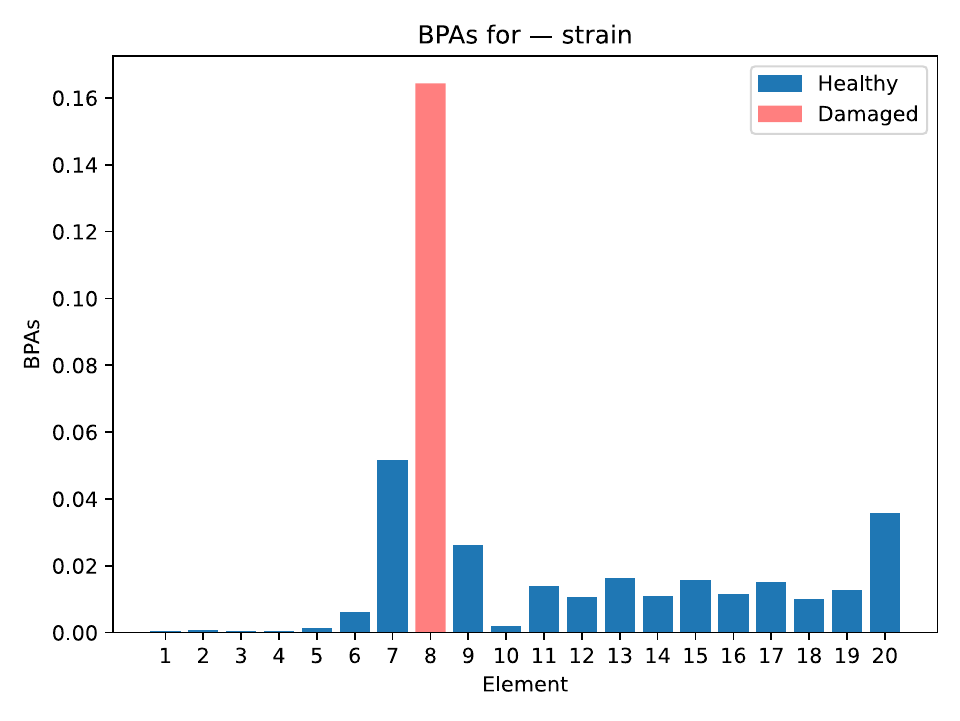%

				\label{fig:sub3}
			\end{subfigure}
			\hfill
			\begin{subfigure}[b]{0.48\textwidth}
	\setlength\figurewidth{\textwidth}%
	\newlength\svgwidth%
	\graphicspath{{figures/DS/svg-inkscape/}}
	\setlength\svgwidth{\figurewidth}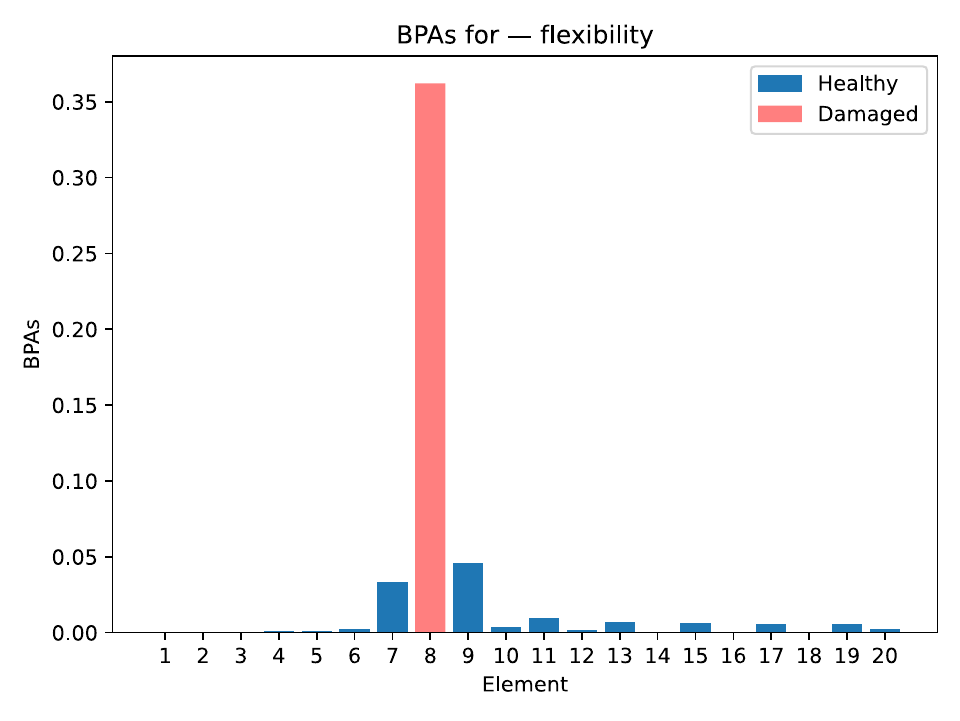%

				\label{fig:sub4}
			\end{subfigure}
			\caption{Overall comparison of the four damage measures}
			\label{fig:ds_sd_all_bpas}
		\end{figure}

	\paragraph{Choosing the best features.}

	As is clear from \Cref{fig:ds_sd_all_bpas}, the flexibiliy- and strain-based measures (MSECR) are the most informative regarding damage localization and hence they are adopted for the fusion using DS evidence theory (\Cref{fig:single_0n_fusion}). Frequency- or curvature-based measures provide little to no useful information for that purpose and hence were not used further.

	\subsubsection{Impact of multiple damage sites on ignorance and fusion quality}
	\label{ssec:multi_damage_ignorance}
	
	When more than one damage site is present, several interrelated effects degrade the precision of DS‐based fusion.  We summarize these below:
	
	\paragraph{i) Increased evidence dispersion.}  
	With $k$ damage sites, the damage index vector $\bm{D}$ becomes more uniform: no single $D_i$ dominates.  From \eqref{eq:entropy}, the normalized entropy $H(\bm{D})$ approaches 1 as the number of comparable nonzero entries grows, so  
	\[
	{conc} = 1 - H(\bm{D}) \to 0.
	\]  
	A low concentration metric means the \emph{distribution uniformity} component $d_i = 1 - {conc}$ rises towards 1 across the board, inflating the ignorance factor $\alpha_i$ as in~\eqref{eq:total_ignorance}.
	
	\paragraph{ii) Higher rank ambiguity.}  
	Multiple similar‐magnitude damage indices induce ties or near‐ties in the ranking step.  Recall $k_i = r_i/(n-1)$ in \Cref{ssec:ignorance_components}.  As more elements share close ranks, the average rank‐based uncertainty $\bar k = \sum_i k_i / n$ increases, again inflating $\alpha_i$.
	
	\paragraph{iii) Feature correlation and conflict.}  
	Strain and flexibility features respond differently to damage location and severity.  With a single crack, their BPA tend to agree (low Dempster conflict).  With multiple cracks, however, feature‐specific BPA often allocate belief to different subsets of $\Theta$, yielding higher conflict mass

	Greater $conflitct$ forces more mass onto the ignorance during normalization, effectively reducing all singleton beliefs $m(\{i\})$.
	
	\paragraph{iv) Dilution of feature confidence.}  
	Each feature’s logistic confidence curve saturates more slowly when multiple moderate $D_i$ values appear.  As $\sigma_i = (1 + \exp(-\lambda_f D_i))^{-1}$ produces mid‐range confidence for several elements, the complement $c_i = 1 - \sigma_i$ remains non‐negligible across many $i$.  Consequently, the \emph{feature confidence} term in \eqref{eq:total_ignorance} contributes more ignorance overall.

	\begin{figure}[tbp]
		\centering
		\begin{subfigure}[b]{0.48\textwidth}
	\setlength\figurewidth{\textwidth}%
	\newlength\svgwidth%
	\graphicspath{{figures/DS/svg-inkscape/}}
	\setlength\svgwidth{\figurewidth}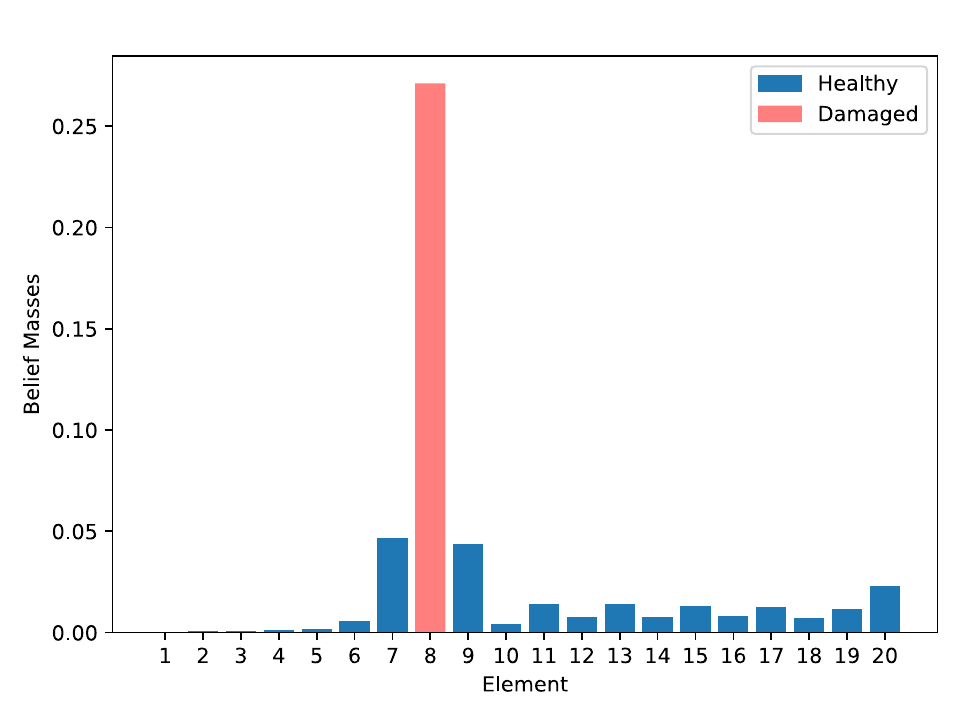%
			\caption{Single damage location}%
			\label{fig:single_0n_fusion}%
		\end{subfigure}
		\hfill%
		\begin{subfigure}[b]{0.48\textwidth}
	\setlength\figurewidth{\textwidth}%
	\newlength\svgwidth%
	\graphicspath{{figures/DS/svg-inkscape/}}
	\setlength\svgwidth{\figurewidth}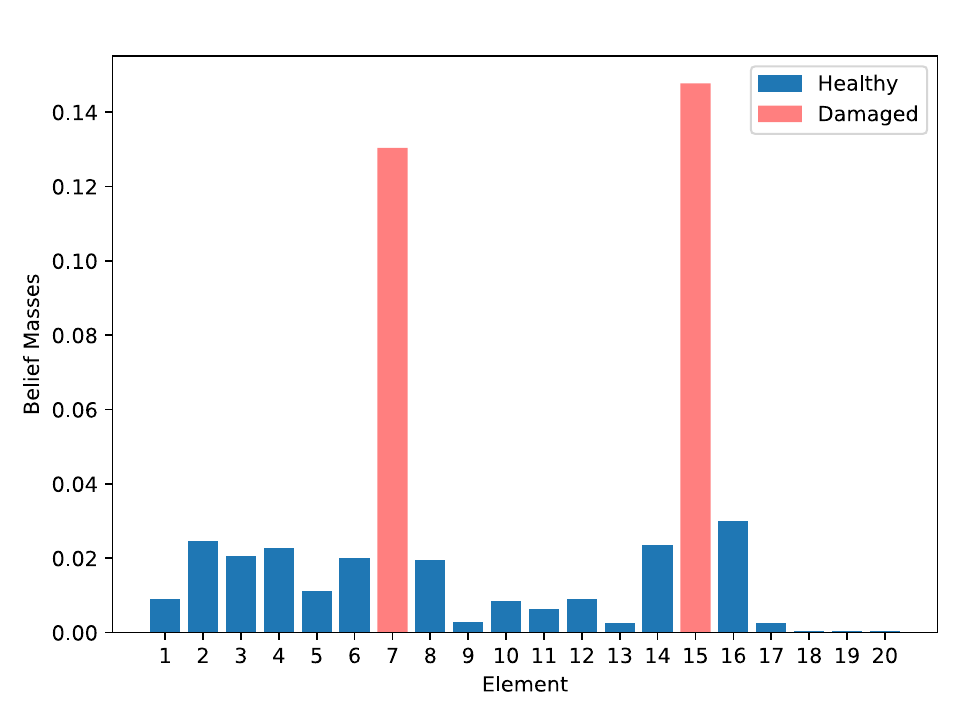%
			\caption{Two damage locations}%
			\label{fig:double_0n_fusion}%
		\end{subfigure}
		\\%
		\begin{subfigure}[b]{0.48\textwidth}
	\setlength\figurewidth{\textwidth}%
	\newlength\svgwidth%
	\graphicspath{{figures/DS/svg-inkscape/}}
	\setlength\svgwidth{\figurewidth}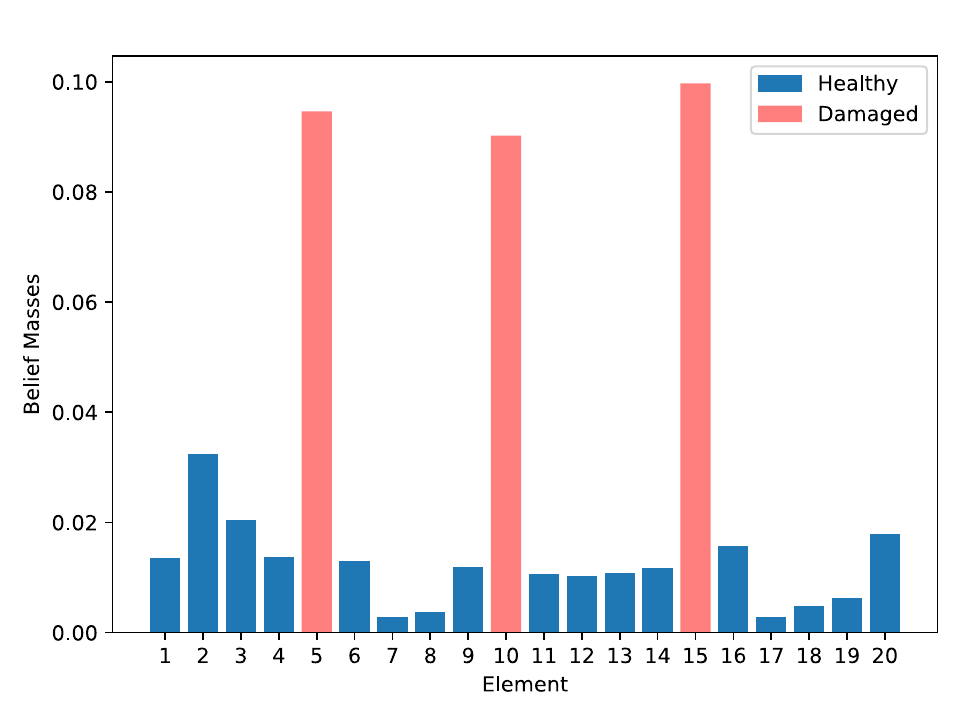%
			\caption{Three damage locations}%
			\label{fig:triple_0n_fusion}%
		\end{subfigure}
		\hfill%
		\begin{subfigure}[b]{0.48\textwidth}
	\setlength\figurewidth{\textwidth}%
	\newlength\svgwidth%
	\graphicspath{{figures/DS/svg-inkscape/}}
	\setlength\svgwidth{\figurewidth}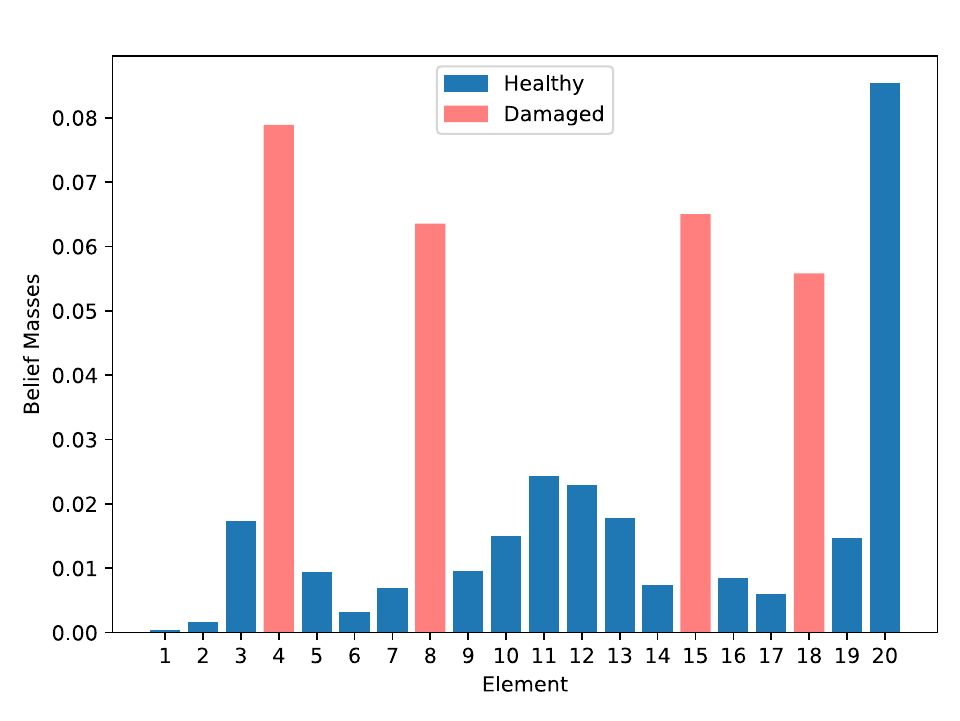%
			\caption{Four damage locations}%
			\label{fig:quad_0n_fusion}%
		\end{subfigure}
		\caption{Effect of multiple damage locations on fused beliefs}
		\label{fig:ds_multiple_damages}
	\end{figure}
	
	All these factors manifest as ‘‘flatter’’ fused belief distributions and lower localization peaks as the number of damage sites grows (cf. \Cref{fig:ds_multiple_damages}). As can be observed when four locations are damaged, the fused belief mass for one healthy element is as high as for the damaged elements and hence can be misleading.
	
	\subsubsection{Effect of measurement noise}
	\label{ssec:ds_noise_effect}

	Measurement noise—particularly in modal eigenvalues and eigenvectors—adversely affects DS‐theoretic fusion by inflating uncertainty and conflict.  In our experiments, we inject zero‐mean Gaussian noise at levels \(\eta \in \{1\%,2\%,3\%,5\%\}\) via
	\begin{align*}
		\tilde\omega_j &= \omega_j\left(1 + \eta\,\mathcal{N}(0,1)\right),\\
		\widetilde\phi_{j} &= \phi_{j} + \eta\,|\phi_{j}|\,\mathcal{N}(0,1),
	\end{align*}
	where \(\omega_j\) and \(\bm{\phi}_{j}(x_i)\) are the \(j\)th eigenfrequency and the mode shape at $x_i$, respectively. While \(\eta=1\%\) reflects a realistic noise floor (cf., e.g., \cite{Grande2007}), we explore up to \(\eta=5\%\) to characterize degradation.

	Noise corrupts the baseline modal features used to compute damage indices \(D_i\).  Random perturbations yield high‐variance estimates \(\hat D_i\), broadening the distribution \(\bm{D}\) and the normalized entropy \(H(\bm{D})\) approaches 1.
	
	All four ignorance components (\(d_i\), \(r'_i\), \(k_i\), \(c_i\)) increase with noise, yielding $\alpha_i$ in \eqref{eq:total_ignorance} that approaches unity as \(\eta\) increases.  Consequently, belief factors \(\beta_i = 1-\alpha_i\) shrink toward zero, and fused belief distributions flatten out.
	
\begin{figure}[tbp]
		\centering
		\begin{subfigure}[b]{0.48\textwidth}
	\setlength\figurewidth{\textwidth}%
	\newlength\svgwidth%
	\graphicspath{{figures/DS/svg-inkscape/}}
	\setlength\svgwidth{\figurewidth}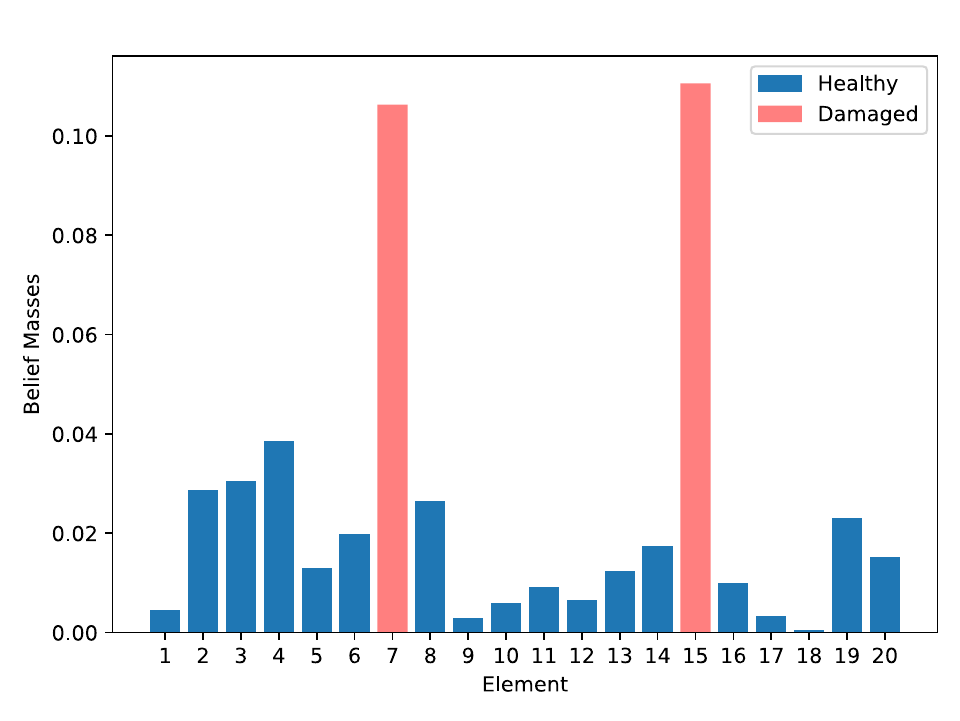%

			\caption{1\% measurement noise}
			\label{fig:doubel_1noise}
		\end{subfigure}
		\hfill
		\begin{subfigure}[b]{0.48\textwidth}
	\setlength\figurewidth{\textwidth}%
	\newlength\svgwidth%
	\graphicspath{{figures/DS/svg-inkscape/}}
	\setlength\svgwidth{\figurewidth}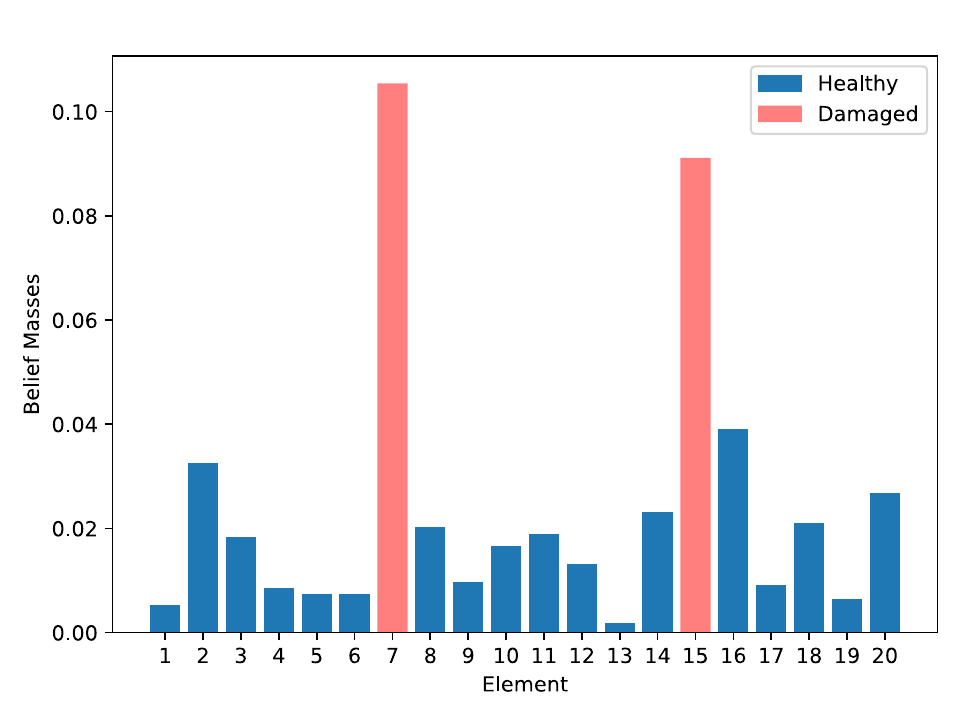%

			\caption{2\% measurement noise}
			\label{fig:doubel_2noise}
		\end{subfigure}
		\vspace{1em}
		\begin{subfigure}[b]{0.48\textwidth}
	\setlength\figurewidth{\textwidth}%
	\newlength\svgwidth%
	\graphicspath{{figures/DS/svg-inkscape/}}
	\setlength\svgwidth{\figurewidth}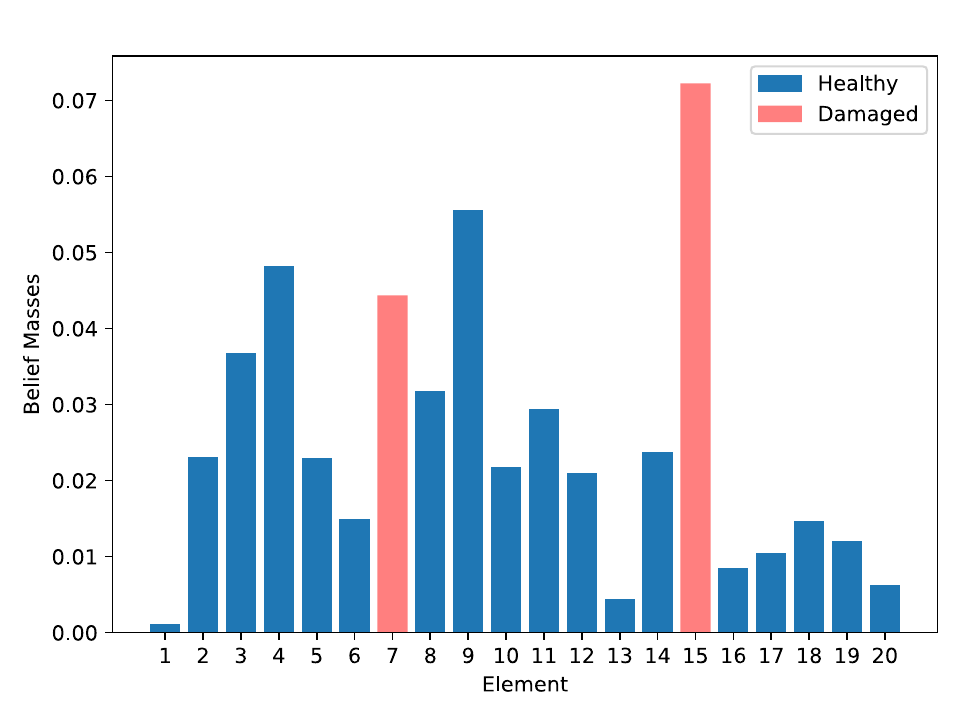%

			\caption{3\% measurement noise}
			\label{fig:doubel_3noise}
		\end{subfigure}
		\hfill
		\begin{subfigure}[b]{0.48\textwidth}
	\setlength\figurewidth{\textwidth}%
	\newlength\svgwidth%
	\graphicspath{{figures/DS/svg-inkscape/}}
	\setlength\svgwidth{\figurewidth}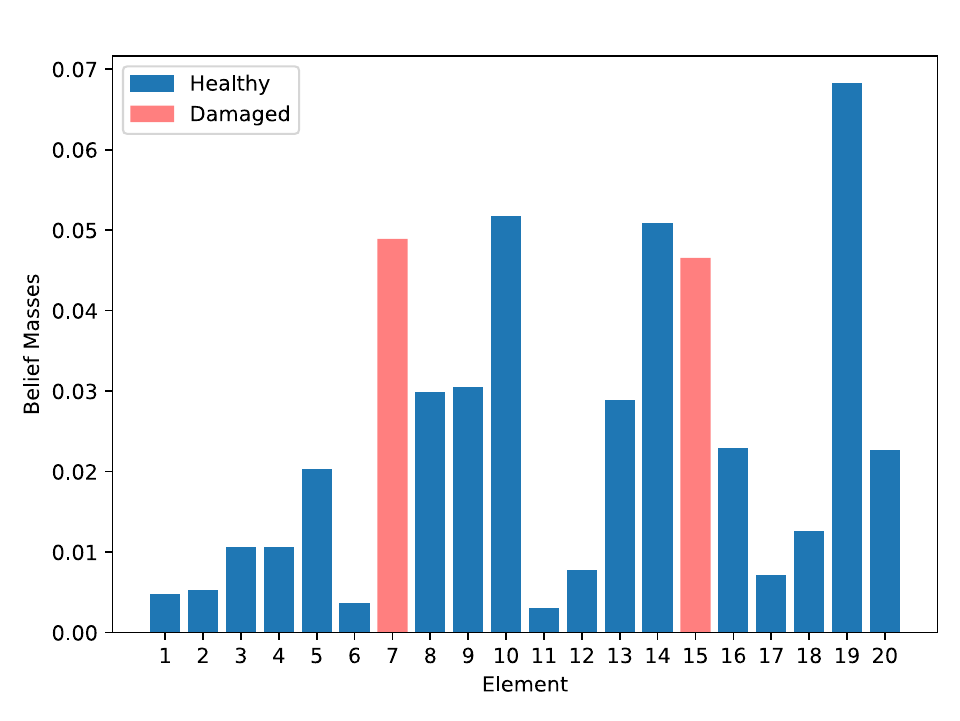%

			\caption{5\% measurement noise}
			\label{fig:doubel_5noise}
		\end{subfigure}
		\caption{Effect of measurement noise on fused damage beliefs}
		\label{fig:ds_noise_effects}
	\end{figure}

	Our four‐case plot series with noise of \(\eta=1,2,3,5\%\) in \Cref{fig:ds_noise_effects} clearly shows progressive deterioration:~At \(\eta=1\%\), localization peaks remain discernible; at \(\eta=5\%\), the damage signal is nearly indistinguishable from noise.  This underscores the need for denoising or robust feature extraction when applying DS fusion in high‐noise environments.

	\subsection{FEM based localization}
	The objective function given by \eqref{eqn:ObjectiveFunctional} in \Cref{ssec:objective_theory} consists of three terms. While the first two terms are more dependent on the global behaviour of the structure, the curvature (third term) is very sensitive to local changes in material parameters. As noted earlier, our inverse Finite Element method formulation minimizes the function using the L-BFGS algorithm. We now present two cases: (i) recovery of a uniform (undamaged) beam, and (ii) identification of localized stiffness reductions. 	
	\paragraph{(i) Healthy beam: convergence to true $E$ values.}
	
\begin{figure}[tbp]
			\centering
			\begin{subfigure}[b]{0.48\textwidth}
	\setlength\figurewidth{\textwidth}%
	\newlength\svgwidth%
	\graphicspath{{figures/iFEM/simple_ifem_ss_healthy_2_terms_obj/svg-inkscape/}}
	\setlength\svgwidth{\figurewidth}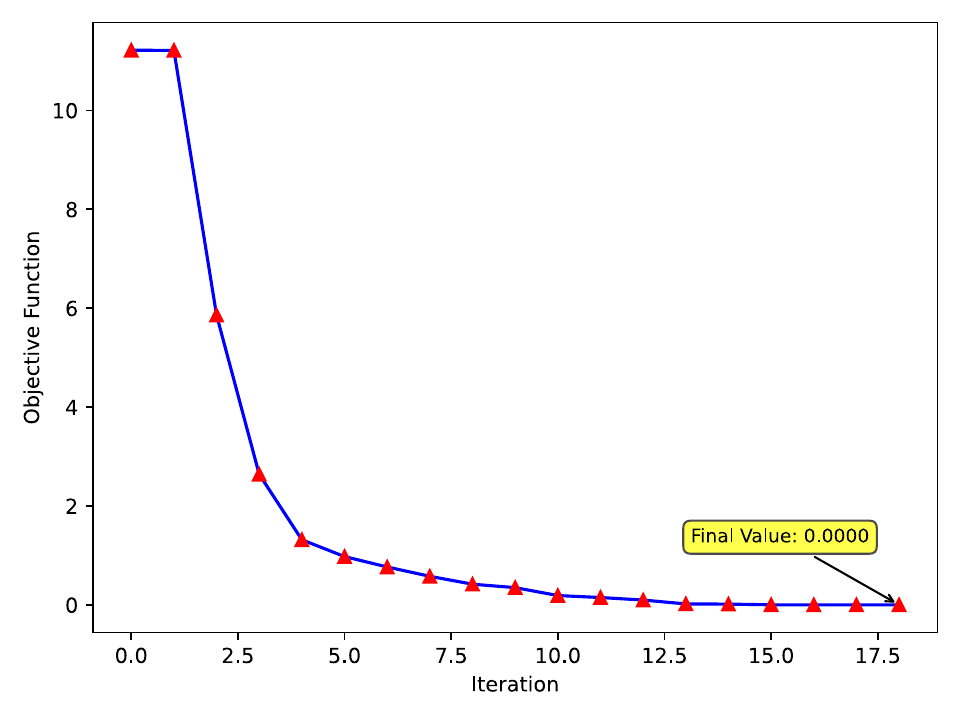%
				\caption{Objective convergence}%
				\label{fig:healthy_ojb}%
			\end{subfigure}
			\hfill
			\begin{subfigure}[b]{0.48\textwidth}
	\setlength\figurewidth{\textwidth}%
	\newlength\svgwidth%
	\graphicspath{{figures/iFEM/simple_ifem_ss_healthy_2_terms_obj/svg-inkscape/}}
	\setlength\svgwidth{\figurewidth}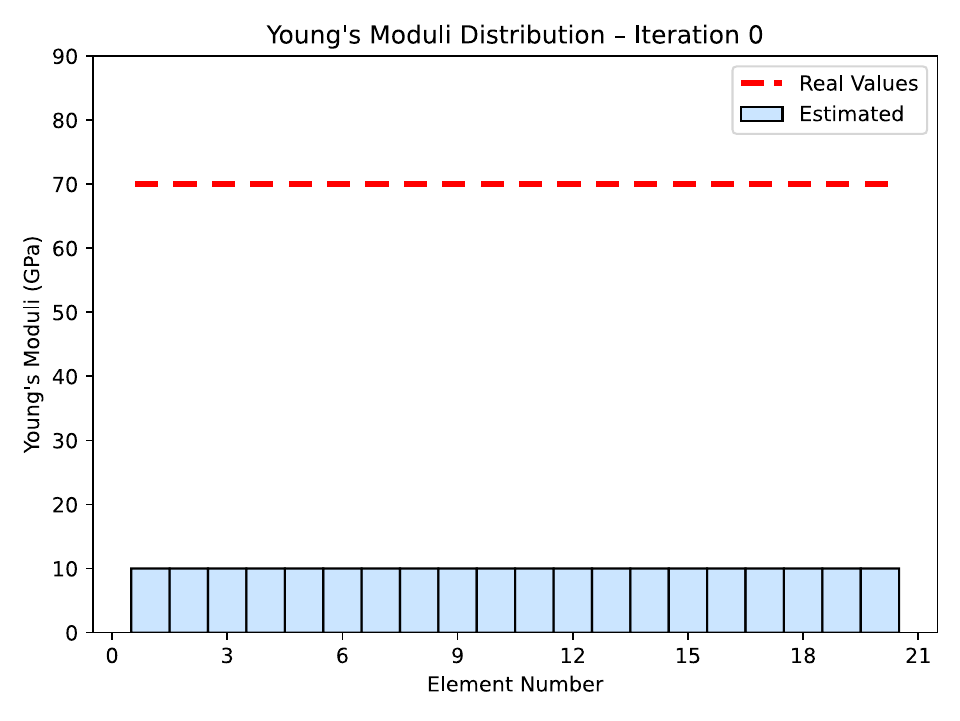%
				\caption{Initial guess }%
				\label{fig:healthy_init_guess}%
			\end{subfigure}
			\vspace{1em}
			\begin{subfigure}[b]{0.48\textwidth}
	\setlength\figurewidth{\textwidth}%
	\newlength\svgwidth%
	\graphicspath{{figures/iFEM/simple_ifem_ss_healthy_2_terms_obj/svg-inkscape/}}
	\setlength\svgwidth{\figurewidth}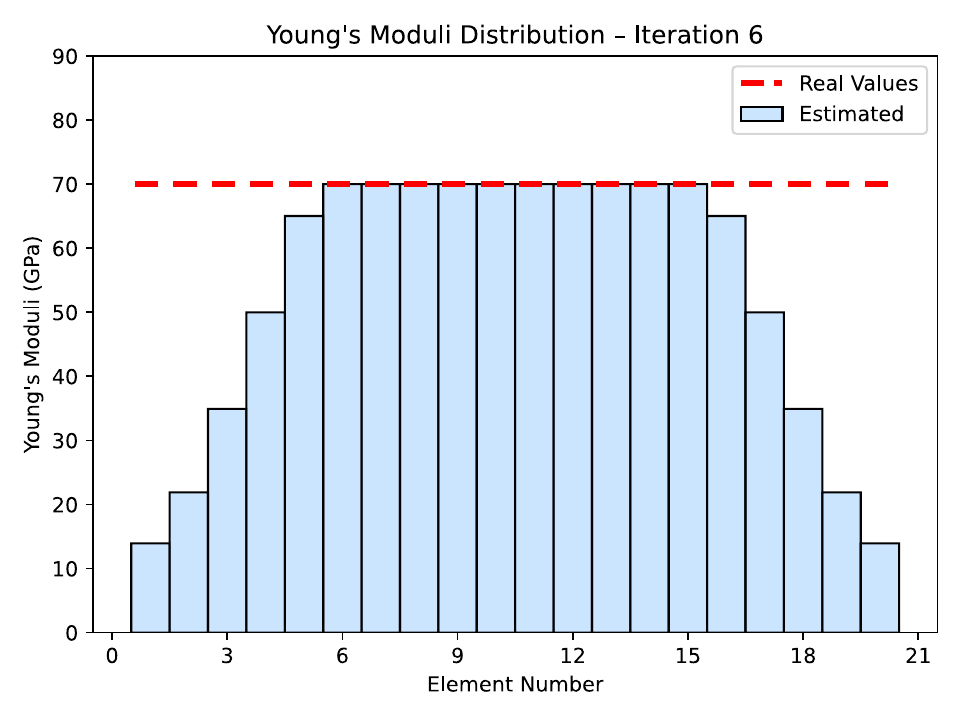%
				\caption{Intermediate distribution}%
				\label{fig:healthy_inter}%
			\end{subfigure}
			\hfill
			\begin{subfigure}[b]{0.48\textwidth}
	\setlength\figurewidth{\textwidth}%
	\newlength\svgwidth%
	\graphicspath{{figures/iFEM/simple_ifem_ss_healthy_2_terms_obj/svg-inkscape/}}
	\setlength\svgwidth{\figurewidth}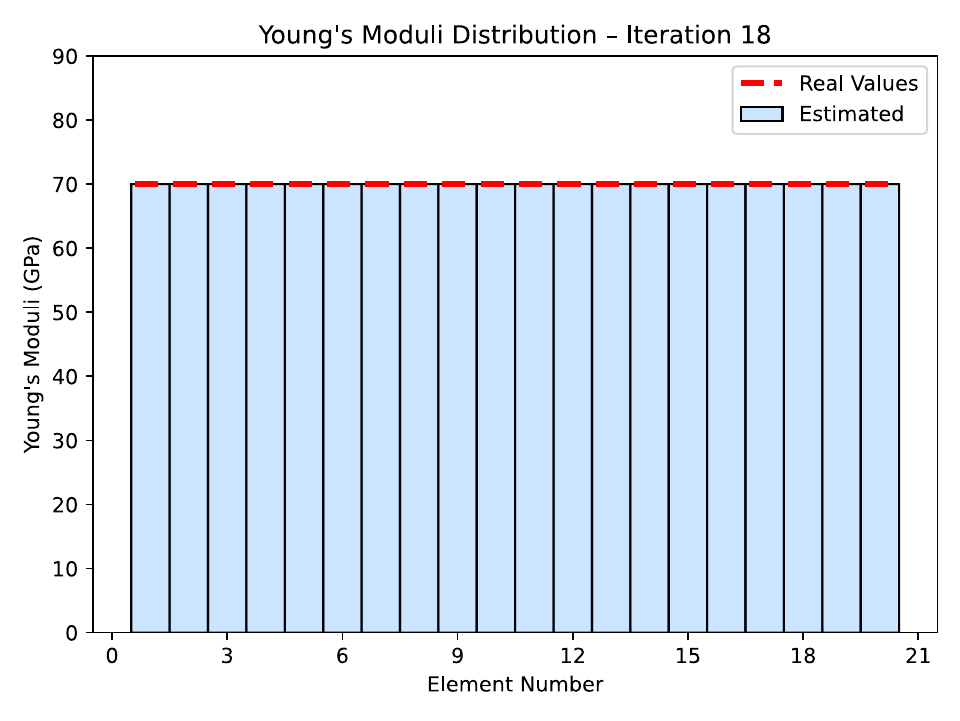%
				\caption{Healthy state}%
				\label{fig:healthy_final}%
			\end{subfigure}
			\caption{FEM based localization convergence to healthy state}
			\label{fig:ifem_healthy}
	\end{figure}

	Starting from a homogeneous initial guess of Young's modulus (e.g., $E_i = 10$\,GPa for all elements), we run a gradient-based optimizer (L-BFGS) to minimize $J(\bm{\theta})$. In the absence of any damage, all elements have the same value of parameters, including their Young's moduli. Thus, the first two terms in \eqref{eqn:ObjectiveFunctional} are sufficient to determine the healthy state as they provide complete information on the global behaviour of the beam. This prevents unnecessary computation of an additional term (curvature) in the objective and therefore also its gradient. 
	\Cref{fig:ifem_healthy} demonstrates that, in the absence of damage, the combined functional provides sufficient global gradient information to guide the search to the correct homogeneous solution.
	
	\paragraph{(ii) Damaged beam: ill-conditioning and local minima.}
	
\begin{figure}[tbp]
			\centering
			\begin{subfigure}[b]{0.48\textwidth}
	\setlength\figurewidth{\textwidth}%
	\newlength\svgwidth%
	\graphicspath{{figures/iFEM/simple_ifem_2_dam_ss/svg-inkscape/}}
	\setlength\svgwidth{\figurewidth}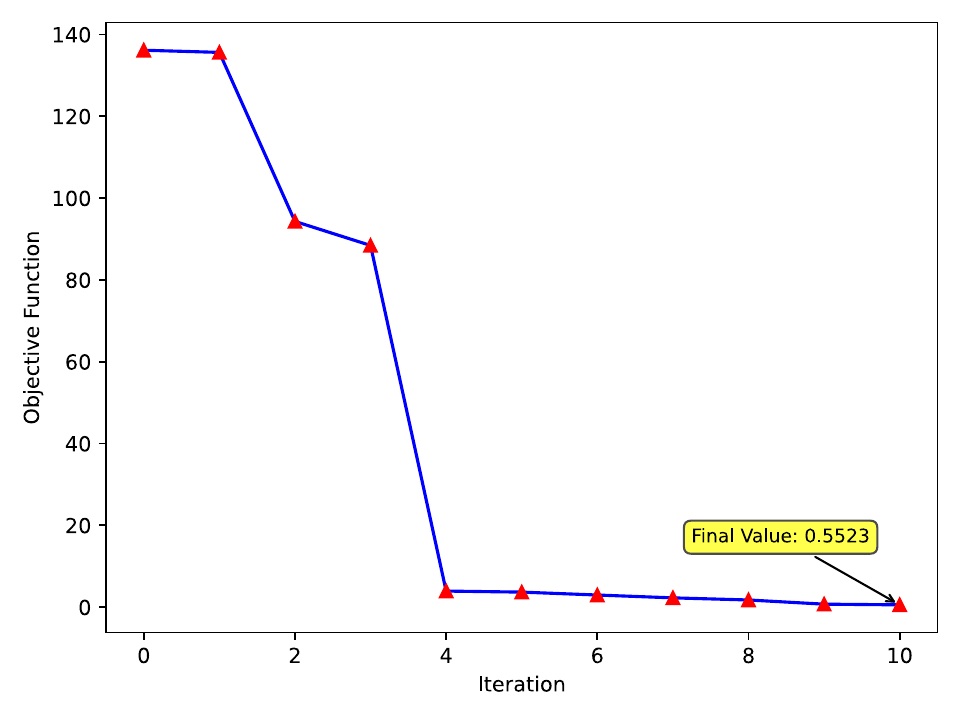%
				\caption{Objective convergence}%
				\label{fig:ifem_objective_damage_2_elem}%
			\end{subfigure}
			\hfill
			\begin{subfigure}[b]{0.48\textwidth}
	\setlength\figurewidth{\textwidth}%
	\newlength\svgwidth%
	\graphicspath{{figures/iFEM/simple_ifem_2_dam_ss/svg-inkscape/}}
	\setlength\svgwidth{\figurewidth}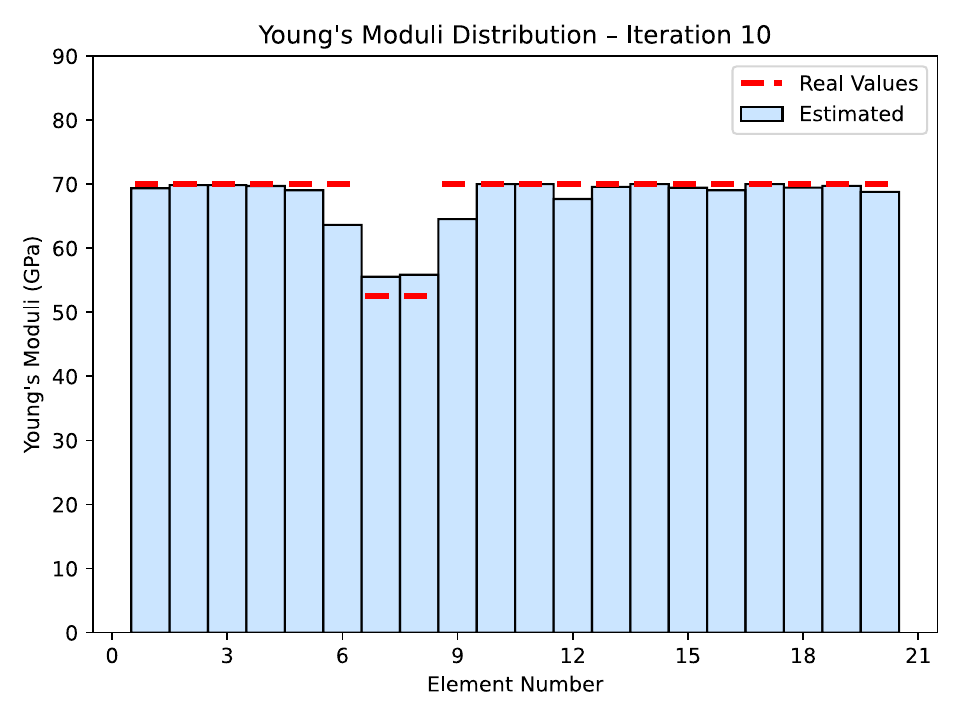%
				\caption{Final damage profile}%
				\label{fig:ifem_youngs_damage_2_elem}%
			\end{subfigure}
			\caption{FEM based convergence for 2 damaged elements}
			\label{fig:ifem_simple_damage}
		\end{figure}

	To simulate damage, we reduce the true Young's moduli $E$ in elements 7 and 8 by 25\% to $E_{damaged} = 52.5$\,GPa, while keeping the rest at healthy values. Using the healthy solution as initial guess\footnote{All subsequent simulations for damage localization are performed using healthy state as the initial guess.}, the optimizer is now tasked with localizing the region of reduced stiffness. Hence all the terms of \eqref{eqn:ObjectiveFunctional} are used, including the regularization to penalize unnecessary changes from the healthy state values. \Cref{fig:ifem_simple_damage} shows the objective convergence and damage profile of the problem. The profile may not match the real values perfectly, but the damaged elements are clearly highlighted but also spread over the adjacent elements. Apart from that there are also small but misleading damages in elements 12 and 20, which of course have no actual damages.

\begin{figure}[tbp]
	 	\centering

	 	\begin{subfigure}[b]{0.48\textwidth}
	\setlength\figurewidth{\textwidth}%
	\newlength\svgwidth%
	\graphicspath{{figures/iFEM/simple_ifem_ss_4_dam/svg-inkscape/}}
	\setlength\svgwidth{\figurewidth}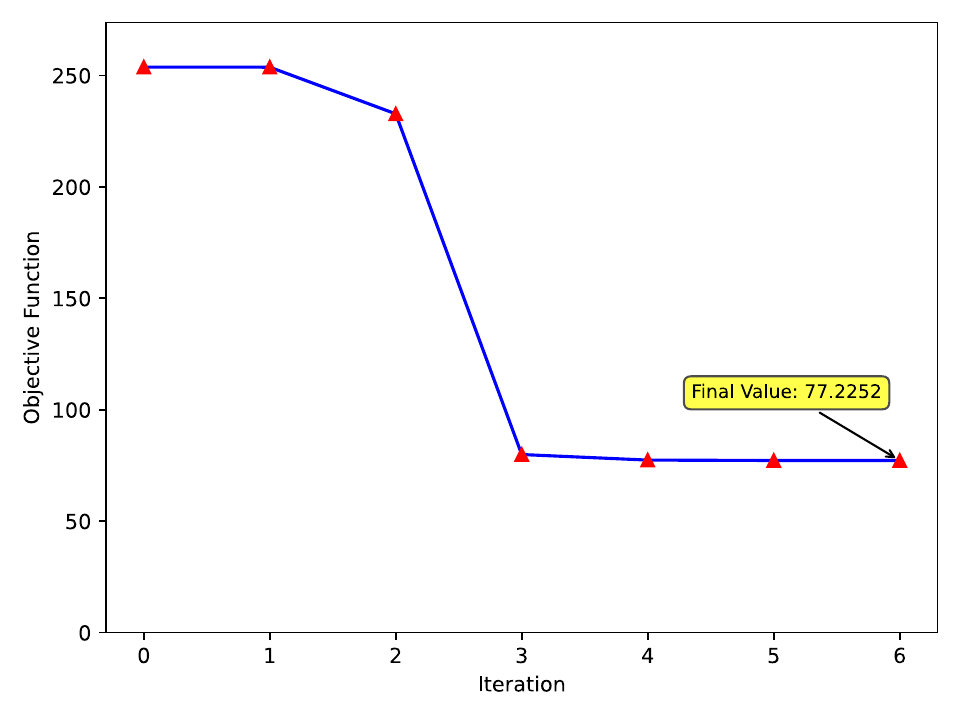%
	 		\caption{Objective convergence}%
	 		\label{fig:ifem_objective_damage_4_elem}%
	 	\end{subfigure}
	 	\hfill
	 	\begin{subfigure}[b]{0.48\textwidth}
	\setlength\figurewidth{\textwidth}%
	\newlength\svgwidth%
	\graphicspath{{figures/iFEM/simple_ifem_ss_4_dam/svg-inkscape/}}
	\setlength\svgwidth{\figurewidth}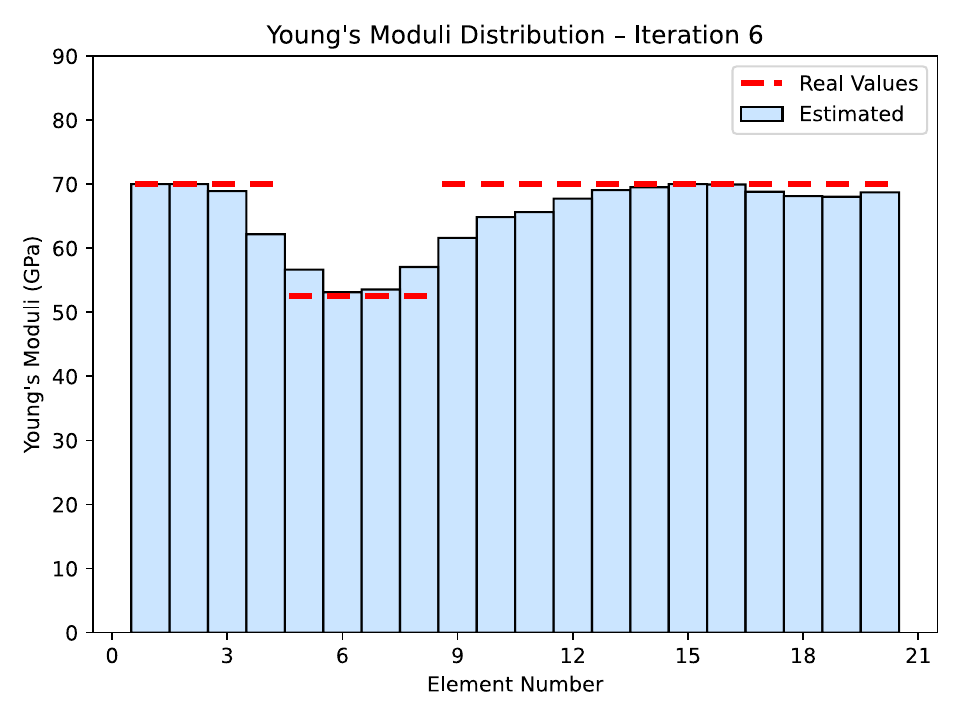%
	 		\caption{Final damage profile}%
	 		\label{fig:ifem_youngs_damage_4_elem}%
	 	\end{subfigure}
	 	\caption{FEM based convergence with 4 damaged elements}
	 	\label{fig:ifem_simple_problem}
	 \end{figure}

	 However, the solution may not always converge smoothly for every possible damage location or any boundary condition, as is clear in \Cref{fig:ifem_simple_problem}.
	 
\begin{figure}[tbp]
		\centering
		\begin{subfigure}[b]{0.48\textwidth}
	\setlength\figurewidth{\textwidth}%
	\newlength\svgwidth%
	\graphicspath{{figures/iFEM/simple_ifem_2_2_dam_lbf/svg-inkscape/}}
	\setlength\svgwidth{\figurewidth}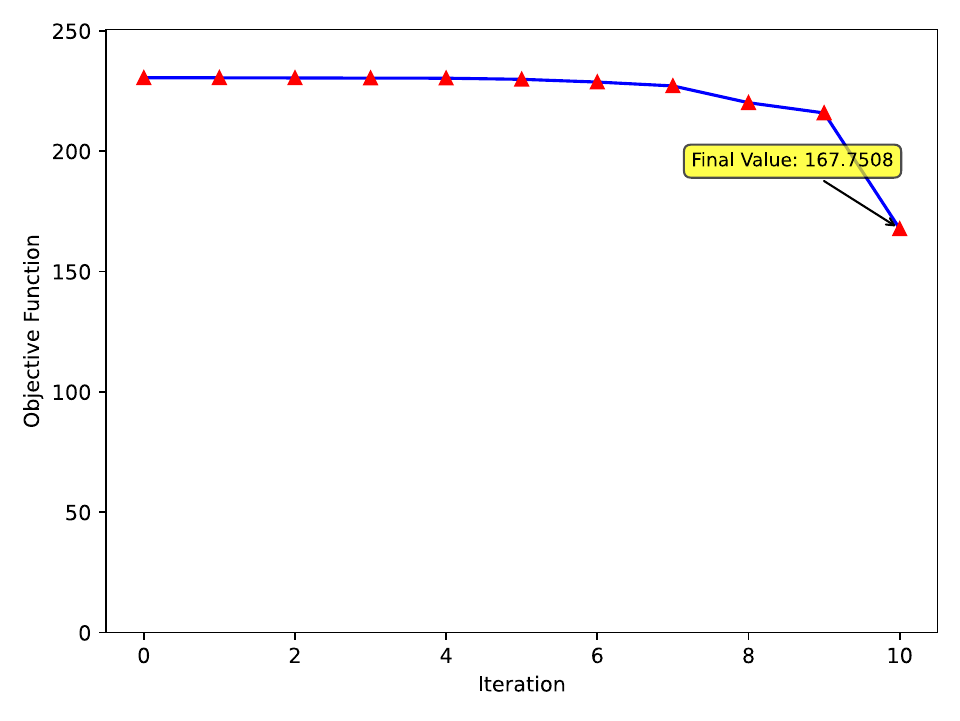%
			\caption{Objective convergence}%
		\end{subfigure}
		\hfill
		\begin{subfigure}[b]{0.48\textwidth}
	\setlength\figurewidth{\textwidth}%
	\newlength\svgwidth%
	\graphicspath{{figures/iFEM/simple_ifem_2_2_dam_lbf/svg-inkscape/}}
	\setlength\svgwidth{\figurewidth}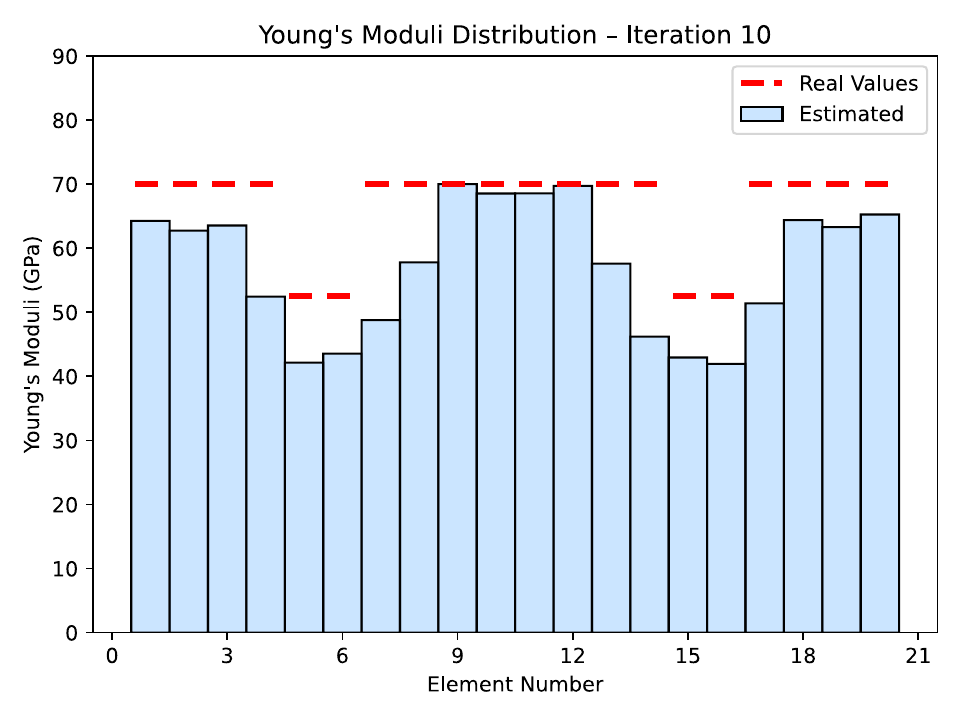%
			\caption{Final damage profile}%
		\end{subfigure}
		\caption{FEM based convergence with 2 separate damage locations}
		\label{fig:simple_ifem_complex}
	\end{figure}

	Although the curvature term of the objective is very sensitive to these local drops, we observe that the objective $J$ can stagnate in a higher-energy ``false'' minimum even after regularization. The results get worse as the damage gets more complex. This is illustrated by \Cref{fig:simple_ifem_complex}, where damage is present at two separate locations. The results only provide very crude information about the damage profile.
	
	 Like many other inverse problems, which are ill-conditioned, non-uniqueness and local minima are a big concern in this case as well. The difficulty in damage localization highlights the need for enhanced optimization strategies.

	\subsubsection{Hierarchical approach}
	
	As explained in \Cref{ssec:hierarchical}, we start with a coarse discretization and gradually refine it. As always, the beam is discretized with 20 elements to keep the comparison consistent. Following notations are used:
	$N_g$ as the initial number of groups, $G_s$ as the initial group size. Also in the convergence plots, the label of the horizontal axis is changed to 'Evaluation number' rather than 'Iteration' to emphasize the fact that the optimizer is restarted at each stage. 
	\paragraph{Case I: $N_g = 5$, $G_s=4$ and damage in one full group.}
	
	This is the simplest case when 5 groups are chosen as design variables in the beginning and one whole group (one cluster) is damaged.
	\begin{figure}[tbp]
		\centering
		\begin{subfigure}[b]{0.48\textwidth}
	\setlength\figurewidth{\textwidth}%
	\newlength\svgwidth%
	\graphicspath{{figures/iFEM/hierarchical/ss_4_dam_lbf/svg-inkscape/}}
	\setlength\svgwidth{\figurewidth}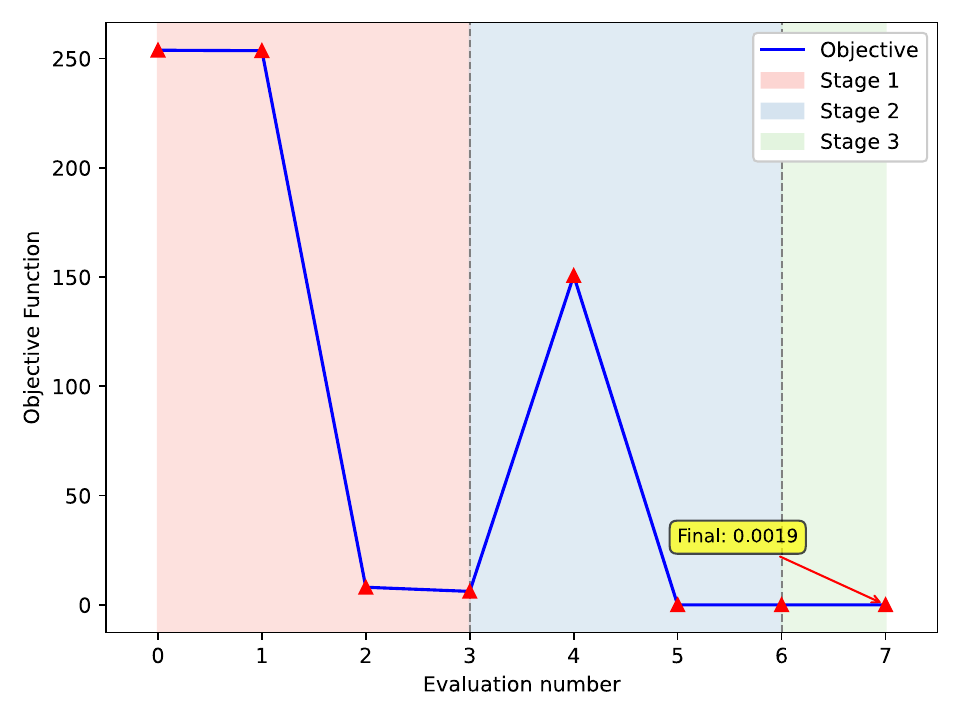%
			\caption{Objective convergence}%
			\label{fig:hier_4_objective_lbf}%
		\end{subfigure}
		\\ %
		\begin{subfigure}[b]{0.48\textwidth}
	\setlength\figurewidth{\textwidth}%
	\newlength\svgwidth%
	\graphicspath{{figures/iFEM/hierarchical/ss_4_dam_lbf/svg-inkscape/}}
	\setlength\svgwidth{\figurewidth}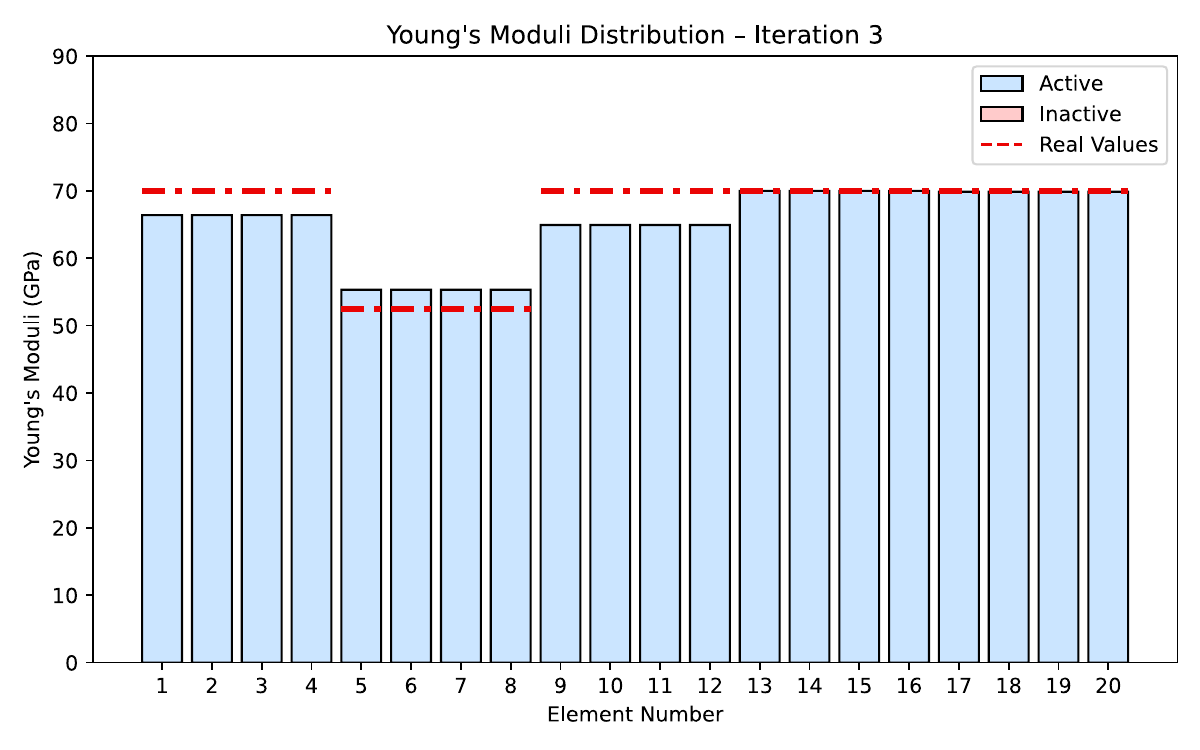%
			\caption{Intermediate result}%
			\label{fig:hier_4_lbf_inter}%
		\end{subfigure}
		\hfill
		\begin{subfigure}[b]{0.48\textwidth}
	\setlength\figurewidth{\textwidth}%
	\newlength\svgwidth%
	\graphicspath{{figures/iFEM/hierarchical/ss_4_dam_lbf/svg-inkscape/}}
	\setlength\svgwidth{\figurewidth}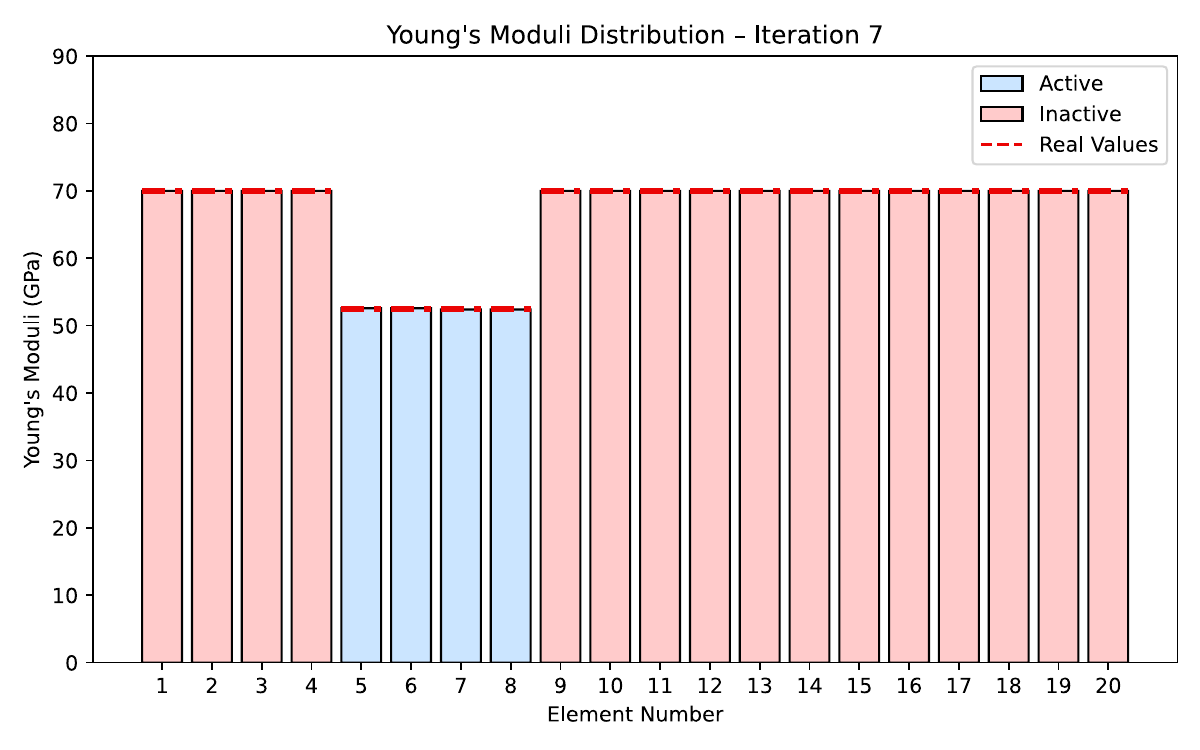%
			\caption{Final damaged profile}%
			\label{fig:hier_4_lbf_final}%
		\end{subfigure}
		\caption{Hierarchical approach for full cluster damage}
		\label{fig:hier_4_lbf}
	\end{figure}
	
	\Cref{fig:hier_4_lbf} shows objective convergence and progression toward the final damaged state. In \Cref{fig:hier_4_objective_lbf}, one can notice the jump in objective value due to stage transfer resulting in freezing of inactive clusters to their original healthy Young's moduli values and division of active clusters. The freezing process often results in a temporary increment (but may also result in a decrement) of the objective value.
	
	\paragraph{Case II: $N_g = 5$, $G_s=4$ and damage in two sub-groups.}
	
\begin{figure}[tbp]
 		\centering
 		\begin{subfigure}[b]{0.48\textwidth}
	\setlength\figurewidth{\textwidth}%
	\newlength\svgwidth%
	\graphicspath{{figures/iFEM/hierarchical/ss_2_2_dam_lbf/eig_und_eq_und_curv/svg-inkscape/}}
	\setlength\svgwidth{\figurewidth}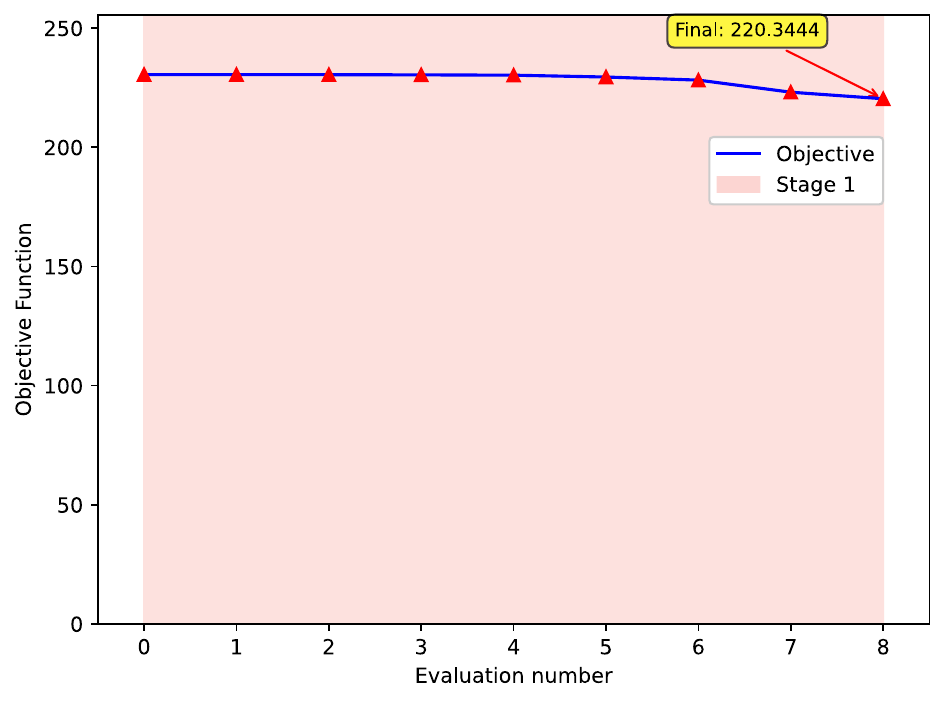%
 			\caption{Objective convergence}%
 			\label{fig:hier_2_full_objective_lbf}%
 		\end{subfigure}
 		\hfill
 		\begin{subfigure}[b]{0.48\textwidth}
	\setlength\figurewidth{\textwidth}%
	\newlength\svgwidth%
	\graphicspath{{figures/iFEM/hierarchical/ss_2_2_dam_lbf/eig_und_eq_und_curv/svg-inkscape/}}
	\setlength\svgwidth{\figurewidth}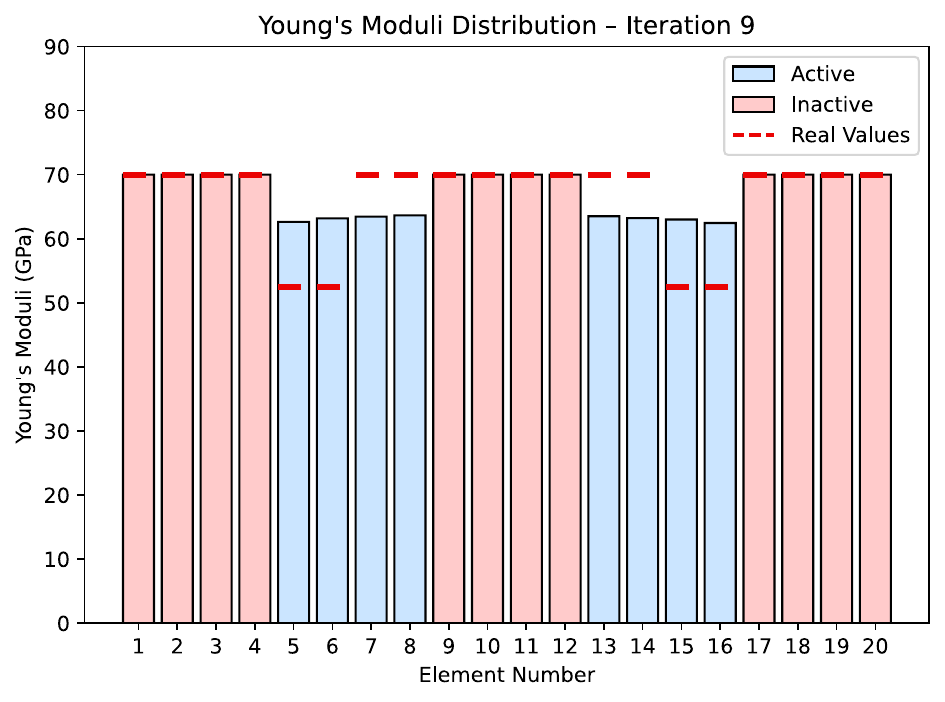%
 			\caption{Final result 1}%
 			\label{fig:hier_2_full_lbf_inter}%
 		\end{subfigure}
 		\caption{Hierarchical L-BFGS approach for sub-cluster damage using full objective}
 		\label{fig:hier_2_full_lbf}
 	\end{figure}

	The hierarchical approach is now tested for a more complex case, where it is not possible to locate the damage in the coarsely discretized state. Instead, the algorithm is required to go for further refinements. Additionally, two locations are damaged. The results for this case can bee seen in \Cref{fig:hier_2_full_lbf} and it is clear that the process has stalled without yielding much information. This is caused by the significant domination and variation in the curvature term in the intermediate results, which prevents the algorithm from taking any further step. This problem has been deliberately chosen to highlight that even the hierarchical approach is not a universal solution of any complex damage localization problem.
	
\begin{figure}[tbp]
		\centering
		\begin{subfigure}[b]{0.48\textwidth}
	\setlength\figurewidth{\textwidth}%
	\newlength\svgwidth%
	\graphicspath{{figures/iFEM/hierarchical/ss_2_2_dam_lbf/eig_und_eq/svg-inkscape/}}
	\setlength\svgwidth{\figurewidth}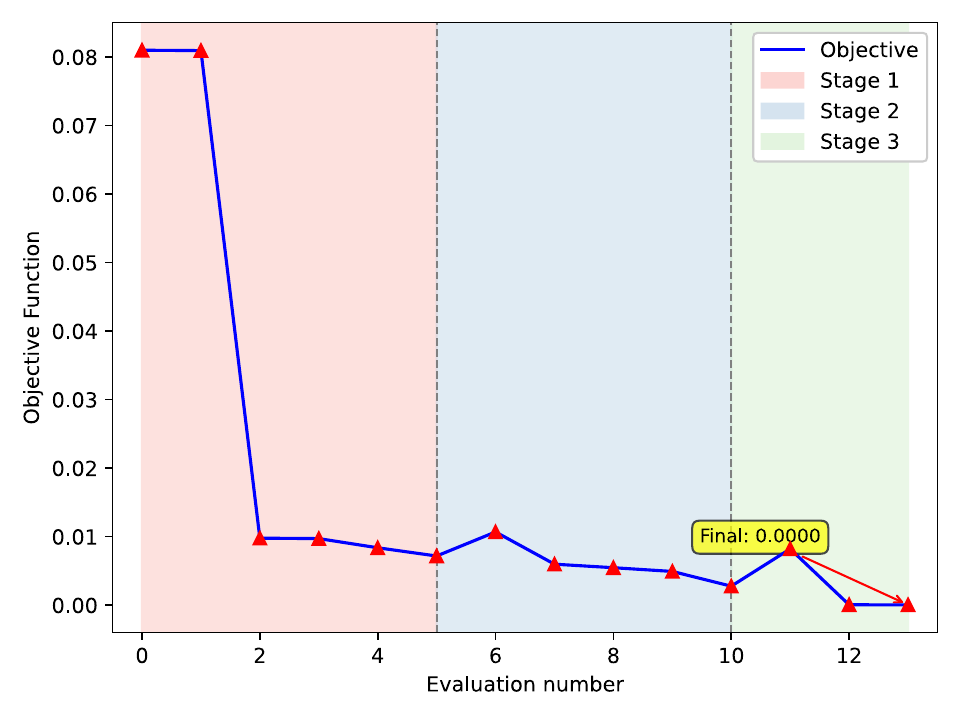%
			\caption{Objective convergence}%
			\label{fig:hier_2_objective_lbf}%
		\end{subfigure}
		\hfill
		\begin{subfigure}[b]{0.48\textwidth}
	\setlength\figurewidth{\textwidth}%
	\newlength\svgwidth%
	\graphicspath{{figures/iFEM/hierarchical/ss_2_2_dam_lbf/eig_und_eq/svg-inkscape/}}
	\setlength\svgwidth{\figurewidth}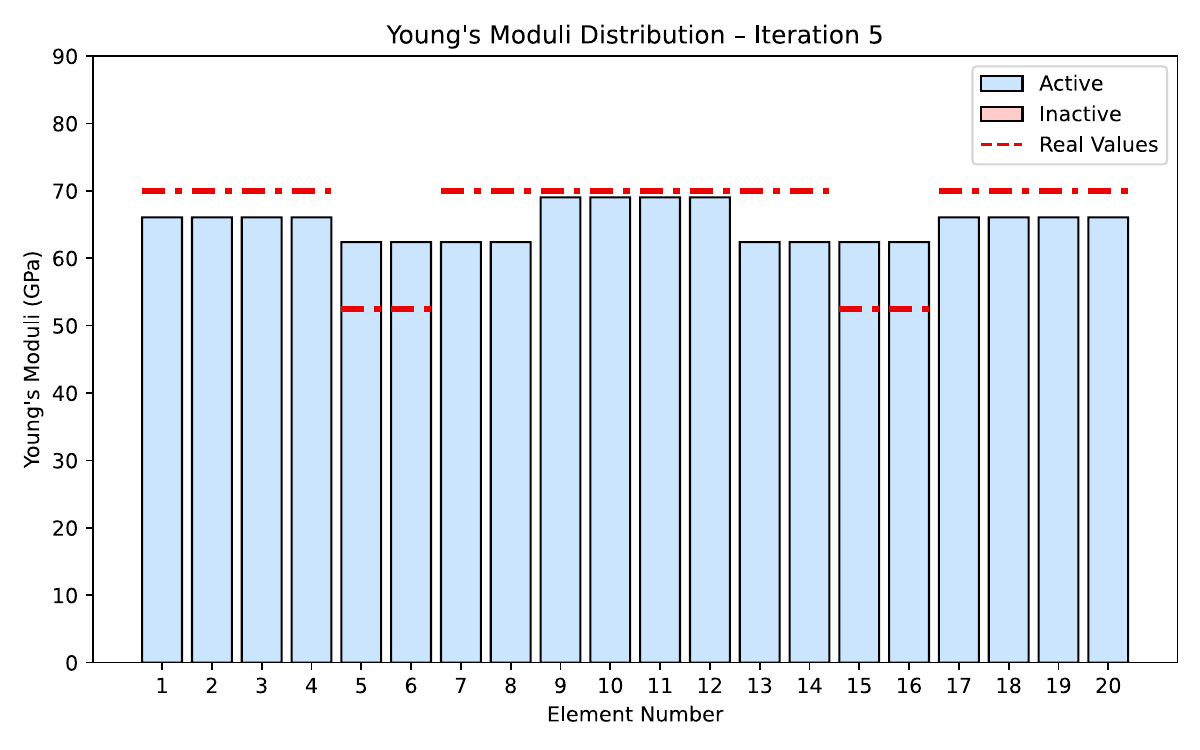%
			\caption{Intermediate result 1}%
			\label{fig:hier_2_lbf_inter}%
		\end{subfigure}
		\\ %
		\begin{subfigure}[b]{0.48\textwidth}
	\setlength\figurewidth{\textwidth}%
	\newlength\svgwidth%
	\graphicspath{{figures/iFEM/hierarchical/ss_2_2_dam_lbf/eig_und_eq/svg-inkscape/}}
	\setlength\svgwidth{\figurewidth}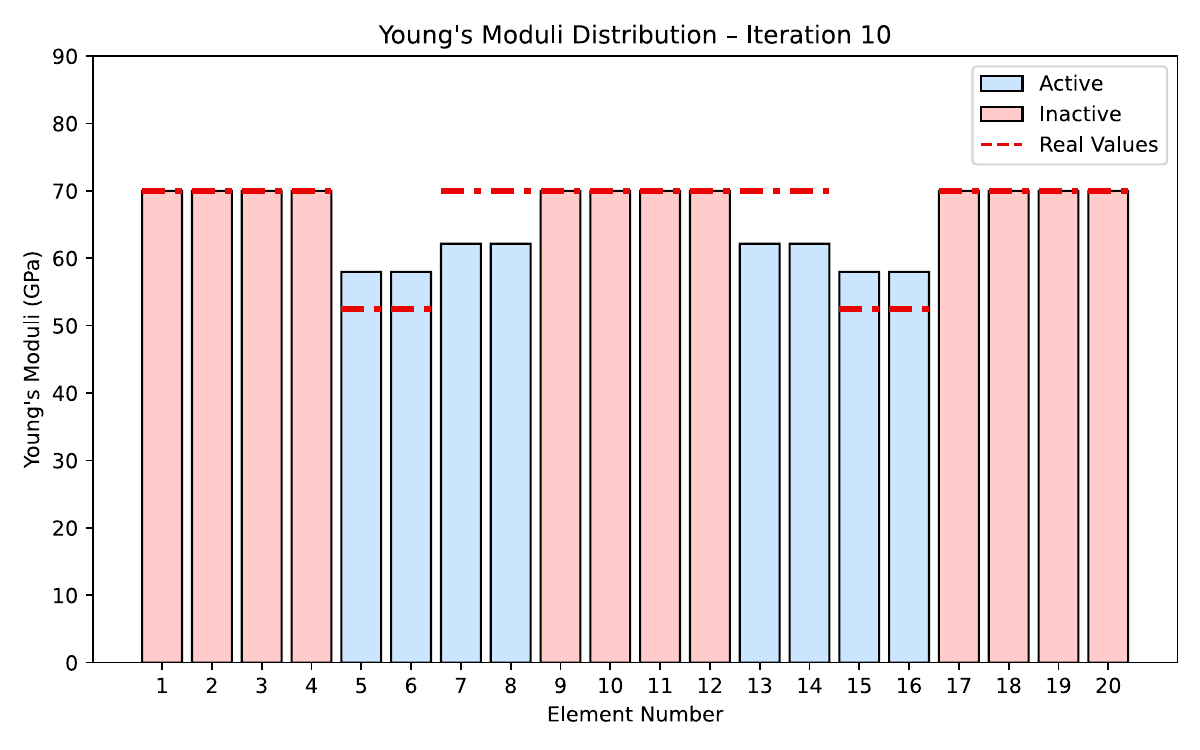%
			\caption{Intermediate result 2}%
			\label{fig:hier_2_lbf_inter_2}%
		\end{subfigure}
		\hfill
		\begin{subfigure}[b]{0.48\textwidth}
	\setlength\figurewidth{\textwidth}%
	\newlength\svgwidth%
	\graphicspath{{figures/iFEM/hierarchical/ss_2_2_dam_lbf/eig_und_eq/svg-inkscape/}}
	\setlength\svgwidth{\figurewidth}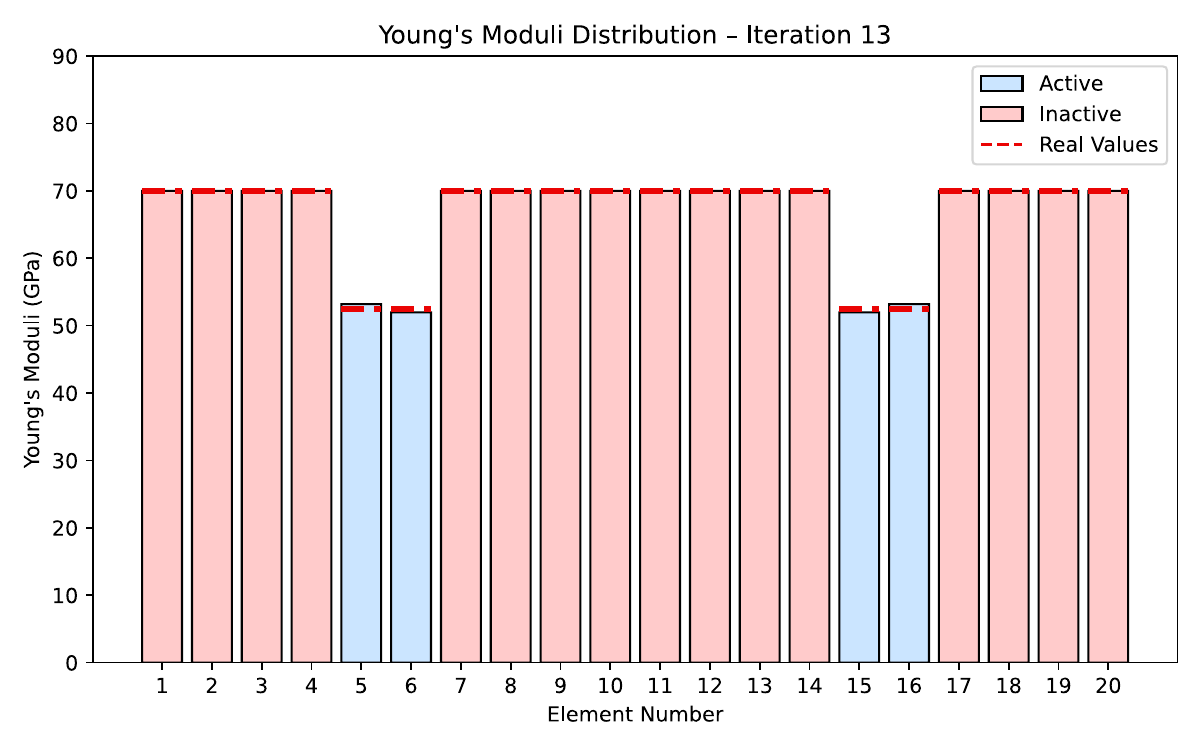%
			\caption{Final damaged profile}%
			\label{fig:hier_2_lbf_final}%
		\end{subfigure}
		\caption{Hierarchical L-BFGS approach for sub-cluster damage omitting curvature error}
		\label{fig:hier_2_lbf}
	\end{figure}

	An interesting result is achieved when the curvature term $E_f(\bm\theta)$ in \eqref{eqn:ObjectiveFunctional} is omitted for this case. One might notice a different scale in the objective (cf. \Cref{fig:hier_2_lbf}). This means that global terms can suffice for convergence when a low-dimensional subproblem instead of a full-dimensional optimization problem is solved. 
	
	As stated earlier, especially in cases like these, the trust-region method outperforms the L-BFGS algorithm without even modifying the original objective function, see \Cref{fig:hier_2_full_lbf}.
	
	\begin{figure}[tbp]
		\centering
		\begin{subfigure}[b]{0.48\textwidth}
	\setlength\figurewidth{\textwidth}%
	\newlength\svgwidth%
	\graphicspath{{figures/iFEM/hierarchical/ss_2_2_dam_trust/eig_und_eq_und_curv/svg-inkscape/}}
	\setlength\svgwidth{\figurewidth}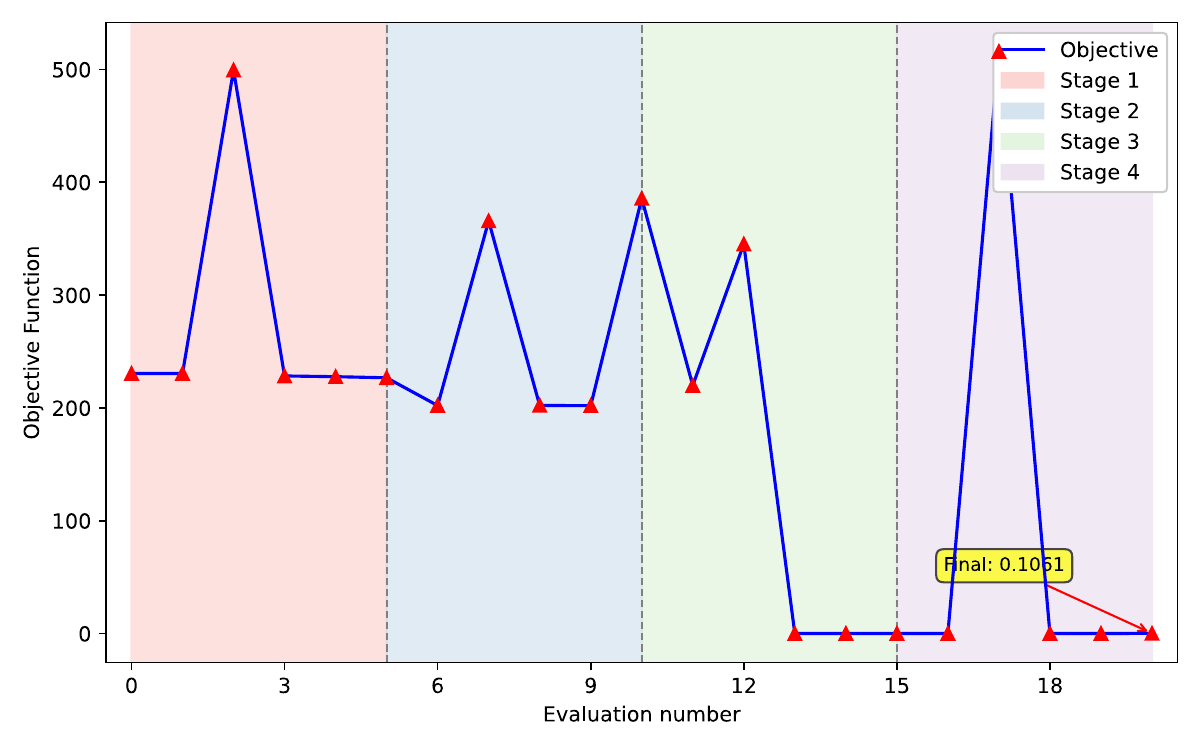%
			\caption{Objective convergence}%
			\label{fig:hier_2_objective_trust}%
		\end{subfigure}
		\hfill
		\begin{subfigure}[b]{0.48\textwidth}
	\setlength\figurewidth{\textwidth}%
	\newlength\svgwidth%
	\graphicspath{{figures/iFEM/hierarchical/ss_2_2_dam_trust/eig_und_eq_und_curv/svg-inkscape/}}
	\setlength\svgwidth{\figurewidth}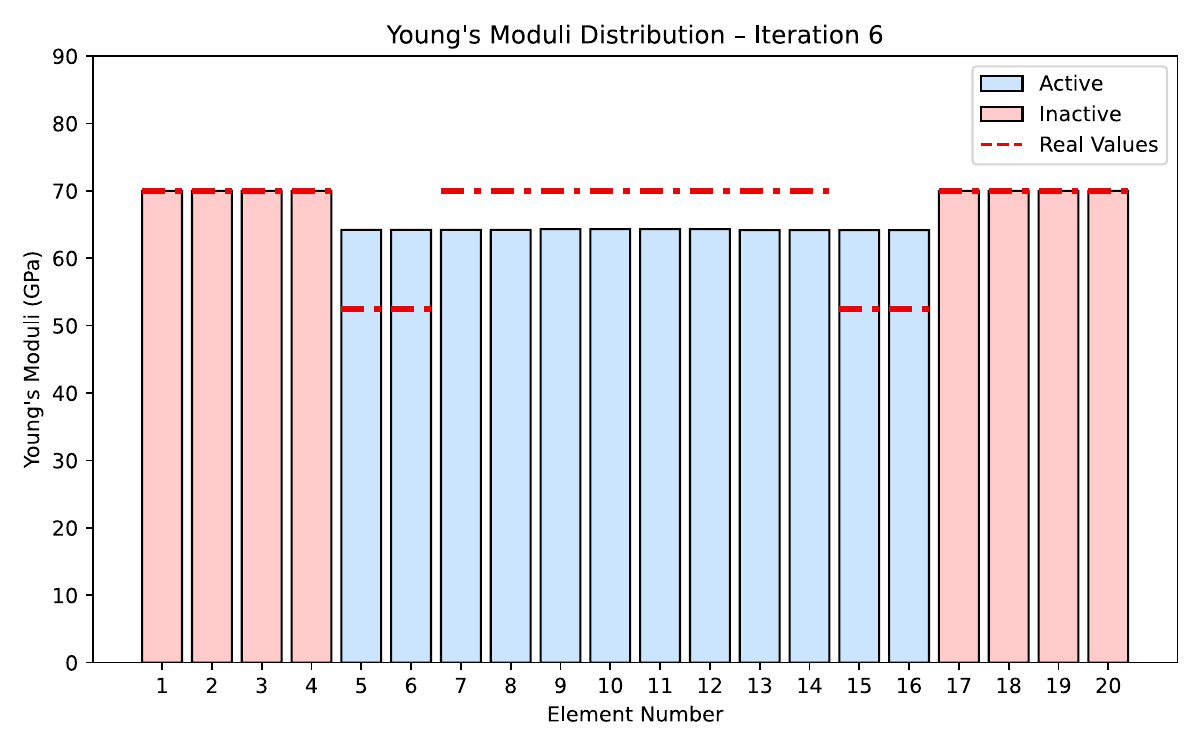%
			\caption{Intermediate result 1}%
			\label{fig:hier_2_trust_inter}%
		\end{subfigure}
		\\ %
		\begin{subfigure}[b]{0.48\textwidth}
	\setlength\figurewidth{\textwidth}%
	\newlength\svgwidth%
	\graphicspath{{figures/iFEM/hierarchical/ss_2_2_dam_trust/eig_und_eq_und_curv/svg-inkscape/}}
	\setlength\svgwidth{\figurewidth}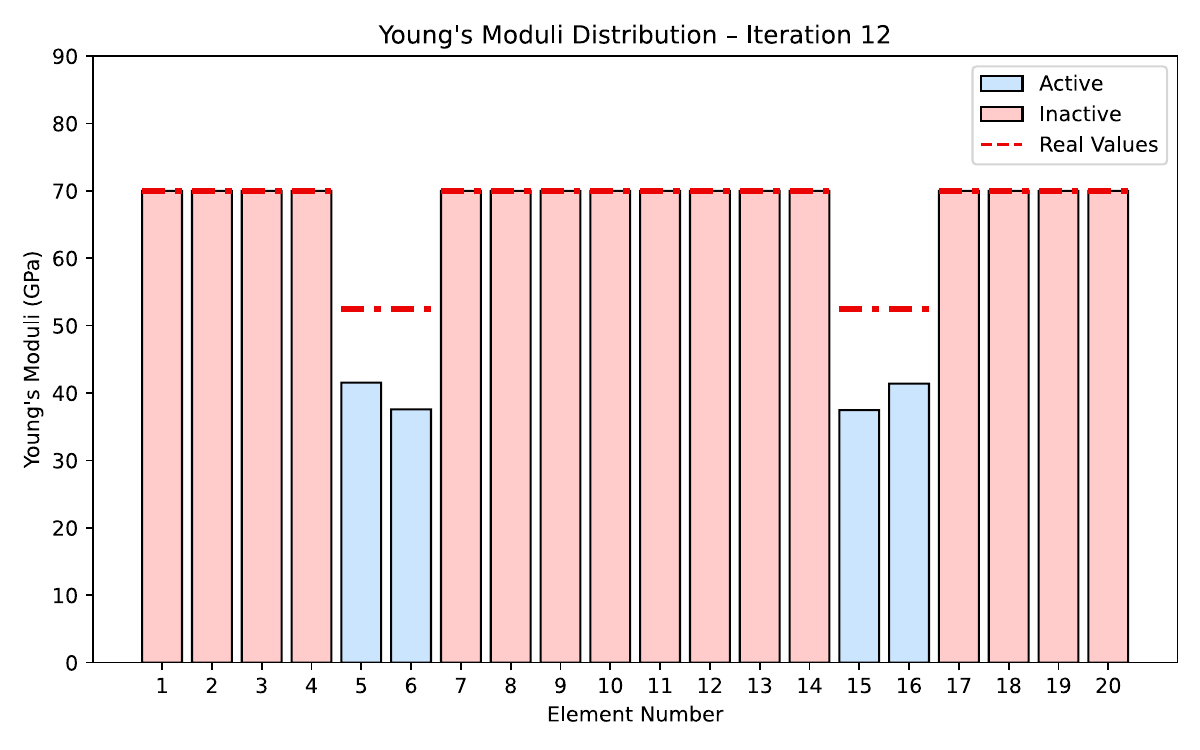%
			\caption{Intermediate result 2}%
			\label{fig:hier_2_trust_inter_2}%
		\end{subfigure}
		\hfill
		\begin{subfigure}[b]{0.48\textwidth}
	\setlength\figurewidth{\textwidth}%
	\newlength\svgwidth%
	\graphicspath{{figures/iFEM/hierarchical/ss_2_2_dam_trust/eig_und_eq_und_curv/svg-inkscape/}}
	\setlength\svgwidth{\figurewidth}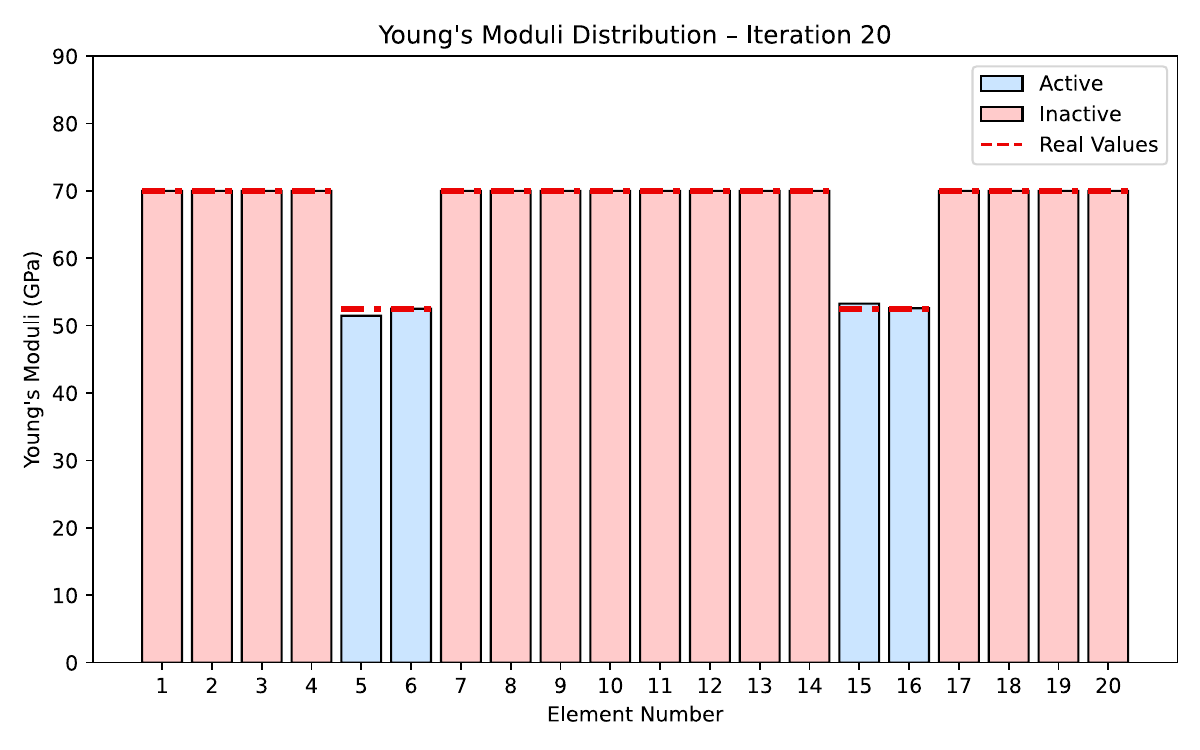%
			\caption{Final damaged profile}%
			\label{fig:hier_2_trust_final}%
		\end{subfigure}
		\caption{Hierarchical trust-region approach for sub-cluster damage using full objective}
		\label{fig:hier_2_trust}
	\end{figure}

	\subsection{Evidence fusion-based candidate filtering and subsequent FEM}
	\label{sec:DS_iFEM_Localization}
	
	One last approach—and by far the most efficient in this study—is to combine Dempster--Shafer (DS) evidence fusion with a focused FEM update on a reduced parameter set. No additional user checks or manual freezing heuristics are required. The single plausibility parameter $\tau$ (taken as a percentage of the highest damage belief in this study), which decides the number of filtered candidates or elements to be passed to the FEM solver, suffices to control model complexity. 
	
	Together, this hybrid DS--FEM framework leverages the strengths of both data-driven evidence fusion and physics-based inversion. It achieves high localization accuracy and reliable quantification with minimal user intervention and computational overhead.

	\begin{figure}[tbp]
		\centering
		\begin{tikzpicture}[
			node distance=10mm,
			every node/.style={font=\footnotesize, align=center},
			startstop/.style={
				rectangle, 
				rounded corners, 
				draw, 
				fill=gray!20, 
				minimum width=3.2cm,
				minimum height=8mm
			},
			process/.style={
				rectangle, 
				draw, 
				fill=blue!10, 
				minimum width=3.2cm,
				minimum height=8mm
			},
			decision/.style={
				diamond, 
				draw, 
				fill=green!10, 
				aspect=1.8,
				inner sep=1pt,
				minimum width=2cm
			},
			flow/.style={-Stealth, thick},
			]
			
			\node (start)     [startstop]                             {Start};
			\node (fusion)    [process, below=of start]               {Run DS evidence fusion};
			\node (filter)    [process, below=of fusion]              {Filter candidates: \\$p_i > \tau$
				\\ or top N candidates};
			\node (ifem)      [process, below=of filter]              {Run FEM on reduced\\parameter set};
			\node (solution)  [startstop, below=of ifem]              {Damage Profile};
			
			\node (tau)      [right=of filter, xshift=5mm]          {$\tau$ (user parameter)};
			\draw [dashed, -Stealth] (tau) -- (filter);
			
			\draw [flow] (start)  -- (fusion);
			\draw [flow] (fusion) -- (filter);
			\draw [flow] (filter) -- (ifem);
			\draw [flow] (ifem)   -- (solution);
			
		\end{tikzpicture}
		\caption{Hybrid DS--FEM damage localization framework.}
		\label{fig:ds_ifem_flowchart}
	\end{figure}
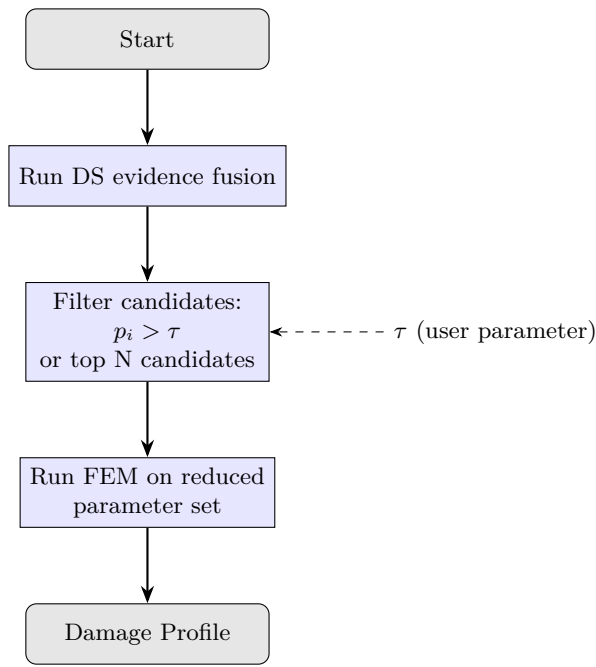

	To increase the reliability and account for measurement noise, the plausibility parameter $\tau$ is taken as 70\% in this study. In order to get a sense of comparison, the same problems discussed previously will now be solved using this method.
	
	\paragraph{Case I: One damage location: four damaged elements.}
	
\begin{figure}[tbp]
		\centering
		\begin{subfigure}[b]{0.48\textwidth}
	\setlength\figurewidth{\textwidth}%
	\newlength\svgwidth%
	\graphicspath{{figures/Fusion/ss_4_lbf_full/svg-inkscape/}}
	\setlength\svgwidth{\figurewidth}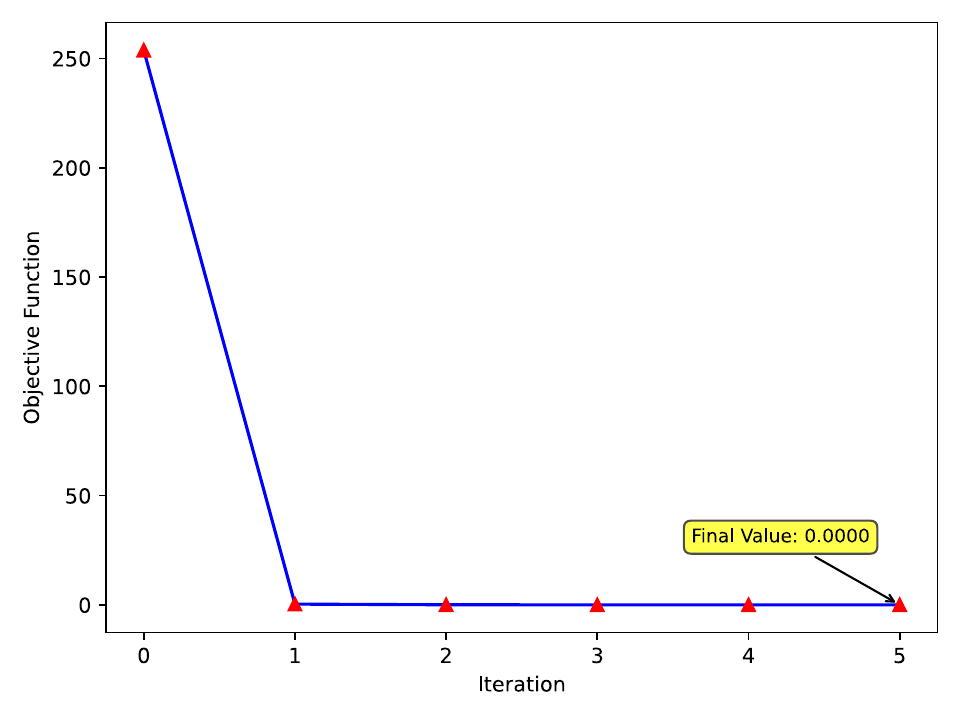%
			\caption{Objective convergence}%
		\end{subfigure}
		\hfill
		\begin{subfigure}[b]{0.48\textwidth}
	\setlength\figurewidth{\textwidth}%
	\newlength\svgwidth%
	\graphicspath{{figures/Fusion/ss_4_lbf_full/svg-inkscape/}}
	\setlength\svgwidth{\figurewidth}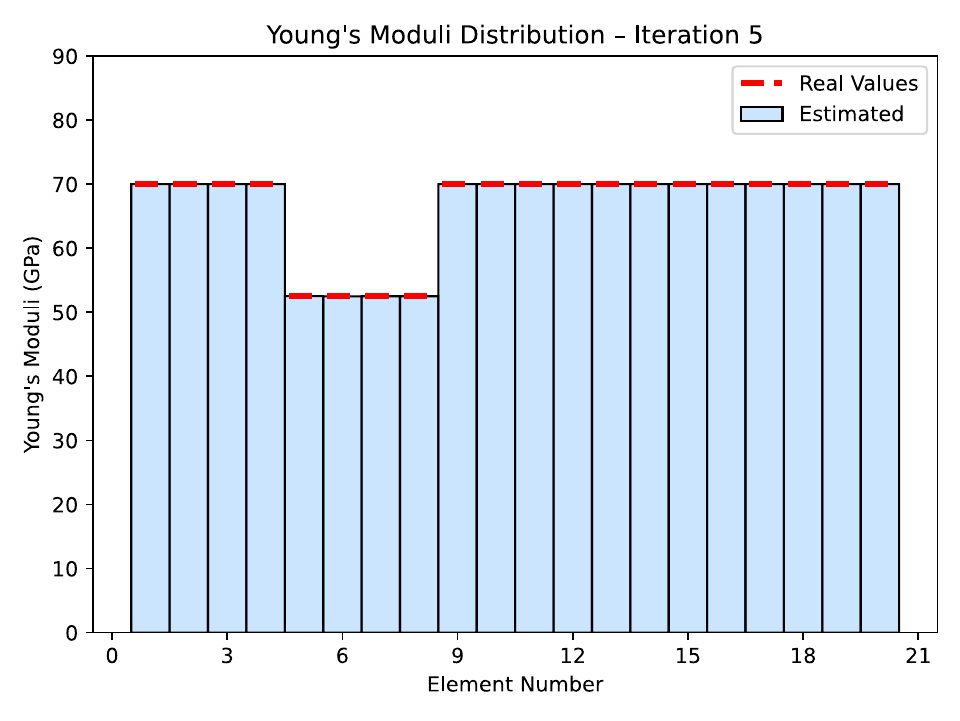%
			\caption{Final damage profile}%
		\end{subfigure}
		\caption{Coupled DS-FEM convergence for case I}
		\label{fig:coupled_4_damage}
	\end{figure}

	This problem has already been solved before using simple FEM based model updation (cf. \Cref{fig:ifem_simple_problem}) and using hierarchical approach (cf. \Cref{fig:ifem_simple_damage}). \Cref{fig:coupled_4_damage} shows the objective and damage profile for the same problem when solved with this coupled approach.
	
	We observe a very low objective value after only one iteration. A solution in accordance with the expected stiffness values is obtained after five iterations, which shows the quick and efficient progress towards detecting present damage.
	
	\paragraph{Case II: Two damage locations: two damaged elements at each location}
	
\begin{figure}[tbp]
		\centering
		\begin{subfigure}[b]{0.48\textwidth}
	\setlength\figurewidth{\textwidth}%
	\newlength\svgwidth%
	\graphicspath{{figures/Fusion/ss_2_2_lbf_full/svg-inkscape/}}
	\setlength\svgwidth{\figurewidth}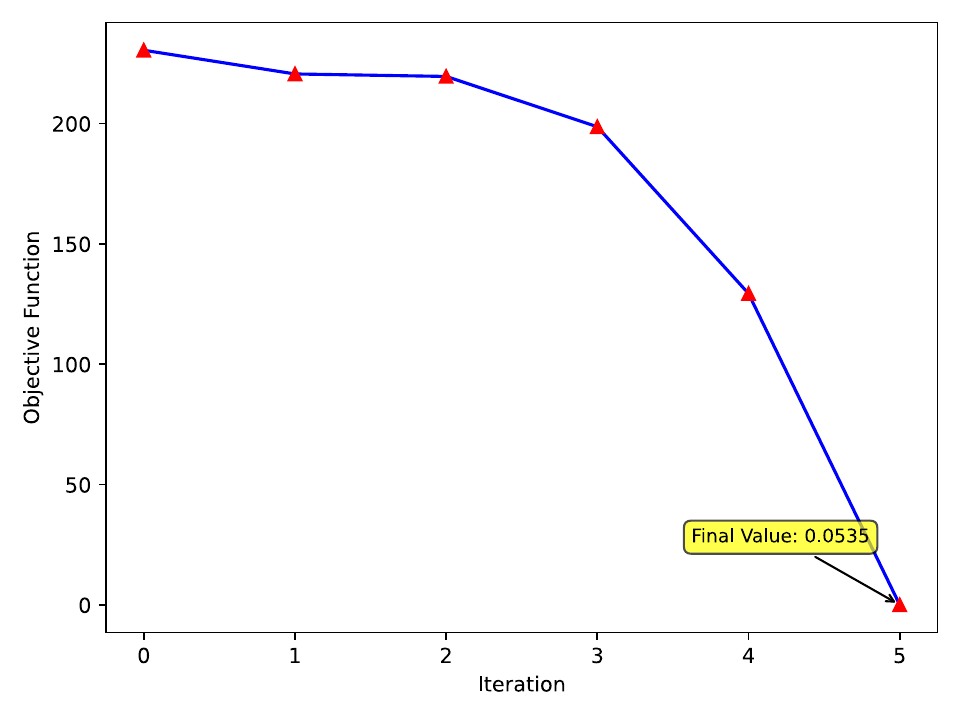%
			\caption{Objective convergence}%
		\end{subfigure}
		\hfill
		\begin{subfigure}[b]{0.48\textwidth}
	\setlength\figurewidth{\textwidth}%
	\newlength\svgwidth%
	\graphicspath{{figures/Fusion/ss_2_2_lbf_full/svg-inkscape/}}
	\setlength\svgwidth{\figurewidth}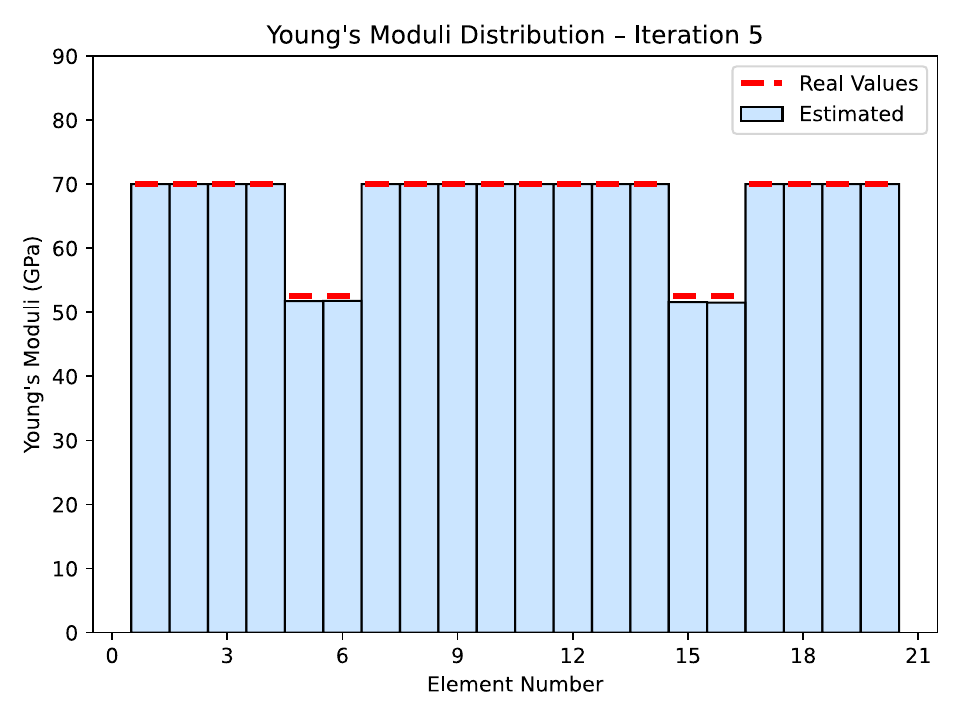%
			\caption{Final damage profile}%
		\end{subfigure}
		\caption{Coupled DS-FEM convergence for case II}
		\label{fig:coupled_2_2_damage}
	\end{figure}

	\Cref{fig:coupled_2_2_damage} shows the objective and damage profile using the coupled DS-FEM approach. The results using simple FEM and the hierarchical approach, have been shown in \Cref{fig:simple_ifem_complex} and \Cref{fig:hier_2_full_lbf}, respectively. 
	
	Again, we observe that a stiffness distribution in accordance with the expected stiffness values is obtained using only five iterations. This shows the superiority of this combined method over the previous methods, which struggle to obtain the correct stiffness distribution.

\section{Conclusions and outlook}
\label{sec:conclusion}
	This study has presented a comprehensive investigation into data- and model-based approaches for structural damage localization and quantification, using a one-dimensional beam as a prototype. The results have demonstrated that while the FEM based model updation methodology is powerful, its straightforward application has been hindered by the inherent ill-conditioning of the inverse problem, which often could not be fully resolved by simple regularization techniques. To overcome this, two novel strategies have been developed: a hierarchical refinement approach and a hybrid technique coupling Dempster-Shafer (DS) evidence theory-based data fusion with FEM.
	
	The hierarchical approach has shown significant improvement over baseline FEM, although its limitations under certain conditions have been identified and documented. In contrast, the hybrid DS-FEM fusion approach has proven superior, achieving the fastest convergence and the most accurate damage profiles.
	
	Subsequent research will involve experimental validation using data from a physical structure to test the framework's performance on real-world damage scenarios.
	
\crefalias{section}{appendix}
\appendix

\section{Calculation of derivatives} \label{app:CalculationDerivatives}

\paragraph{Derivative of frequency-shift error.}
By using $\lambda_j^{mod}(\bm{\theta}) = \left( \omega_j^{mod} (\bm{\theta}) \right)^2$ we directly obtain
	\begin{align*}
		\begin{aligned}
			\frac{\partial E_f(\bm{\theta})}{\partial \theta_{k}}
			&= \sum_{j=1}^{m} 2 \cdot \frac{\omega_{j}^{exp}-\omega_{j}^{mod}(\bm{\theta})}{\left( \omega_{j}^{exp} \right)^2} \cdot \left(- \frac{\partial \omega_{j}^{mod}(\bm{\theta})}{\partial \theta_{k}} \right) \\
			&= -2 \sum_{j=1}^{m} \frac{\omega_{j}^{exp}-\omega_{j}^{mod}(\bm{\theta})}{\left( \omega_{j}^{exp} \right)^2} \cdot \frac{\partial \sqrt{\lambda_{j}^{mod}(\bm{\theta})}}{\partial \theta_{k}} \\
			&= -2 \sum_{j=1}^{m} \frac{\omega_{j}^{exp}-\omega_{j}^{mod}(\bm{\theta})}{\left( \omega_{j}^{exp} \right)^2 } \cdot \frac{1}{2 \sqrt{\lambda_{j}^{mod}(\bm{\theta})} } \frac{\partial \lambda_{j}^{mod}(\bm{\theta})}{\partial \theta_{k}} \\
			&= - \sum_{j=1}^{m} \frac{\omega_{j}^{exp}-\omega_{j}^{mod}(\bm{\theta})}{\left( \omega_{j}^{exp} \right)^2 \omega_{j}^{mod}(\bm{\theta}) } \cdot \frac{\partial \lambda_{j}^{mod}(\bm{\theta})}{\partial \theta_{k}}.
		\end{aligned}
\end{align*}

\paragraph{Derivative of the governing-equation residual.}
As the model mass matrix is independent of $\bm{\theta}$ in our case, the derivative reads
	\begin{align*}
		\begin{aligned}
			&\frac{\partial E_g(\bm{\theta})}{\partial \theta_{k}}
			= \sum_{j=1}^{m} 2 \cdot \left(1 - 
			\frac{\left( \bm{\phi}_{j}^{exp}\right)^\top \bm{K}^{mod}(\bm{\theta}) \,\bm{\phi}_{j}^{exp}}
			{\left( \omega_{j}^{exp} \right)^{2}\,\left( \bm{\phi}_{j}^{exp}\right)^\top \bm{M}^{mod} \,\bm{\phi}_{j}^{exp}}
			\right) \\
			&\hphantom{\frac{\partial E_g(\bm{\theta})}{\partial \theta_{k}}=\sum_{j=1}^{m} 2 \cdot}\cdot \frac{\partial}{\partial \theta_{k}} \left(1 - 
			\frac{\left( \bm{\phi}_{j}^{exp}\right)^\top \bm{K}^{mod}(\bm{\theta}) \,\bm{\phi}_{j}^{exp}}
			{\left( \omega_{j}^{exp} \right)^{2}\,\left( \bm{\phi}_{j}^{exp}\right)^\top \bm{M}^{mod} \,\bm{\phi}_{j}^{exp}}
			\right) \\
			&= 2 \sum_{j=1}^{m} \left(1 - 
			\frac{\left( \bm{\phi}_{j}^{exp}\right)^\top \bm{K}^{mod}(\bm{\theta}) \,\bm{\phi}_{j}^{exp}}
			{\left( \omega_{j}^{exp} \right)^{2}\,\left( \bm{\phi}_{j}^{exp}\right)^\top \bm{M}^{mod} \,\bm{\phi}_{j}^{exp}}
			\right) \cdot \left( - 
			\frac{\left( \bm{\phi}_{j}^{exp}\right)^\top \frac{\partial \bm{K}^{mod}(\bm{\theta})}{\partial \theta_{k}} \,\bm{\phi}_{j}^{exp}}
			{\left( \omega_{j}^{exp} \right)^{2}\,\left( \bm{\phi}_{j}^{exp}\right)^\top \bm{M}^{mod} \,\bm{\phi}_{j}^{exp}} \right).
		\end{aligned}
\end{align*}

\paragraph{Derivative of the curvature-error term.}
The derivative is given as
	\begin{align*}
		\begin{aligned}
			\frac{\partial E_c(\bm{\theta})}{\partial \theta_{k}} 
			&= \sum_{j=1}^{m} \sum_{i=1}^{N} 2 \cdot \frac{\kappa_{j}^{exp}(x_i)-\kappa_{j}^{mod}(x_i,\bm{\theta})}
			{\max\left( \left| \kappa_{j}^{exp}(x_i) \right|, \epsilon \right)^2} \cdot \left( - \frac{\partial \kappa_{j}^{mod}(x_i,\bm{\theta})}{\partial \theta_{k}} \right) \\
			&= -2 \sum_{j=1}^{m} \sum_{i=1}^{N} \frac{\kappa_{j}^{exp}(x_i)-\kappa_{j}^{mod}(x_i,\bm{\theta})}
			{\max\left( \left| \kappa_{j}^{exp}(x_i) \right|, \epsilon \right)^2} \cdot \frac{\partial}{\partial \theta_{k}} \left( \left. \frac{\partial^2}{{\partial x}^2} \left( \bm{\phi}_{j}^{mod}(x,\bm{\theta}) \right) \right|_{x=x_i} \right).
		\end{aligned}
	\end{align*}
	As the eigenvector~$\bm{\phi}_{j}^{mod}(x,\bm{\theta})$ is obtained numerically using FEM, we obtain a discretized version of the eigenvector. Assumung sufficient continuity, the partial derivatives {w.r.t.}~$x$ and $\theta_k$ can be exchanged, which gives
	\begin{align*}
		\frac{\partial E_c(\bm{\theta})}{\partial \theta_{k}}  = -2 \sum_{j=1}^{m} \sum_{i=1}^{N} \frac{\kappa_{j}^{exp}(x_i)-\kappa_{j}^{mod}(x_i,\bm{\theta})}
		{\max\left( \left| \kappa_{j}^{exp}(x_i) \right|, \epsilon \right)^2} \cdot \left. \frac{\partial^2}{{\partial x}^2} \left( \frac{\partial}{\partial \theta_{k}} \left( \bm{\phi}_{j}^{mod}(x,\bm{\theta}) \right) \right) \right|_{x=x_i}.
	\end{align*}
	The second derivative of $\bm{\phi}_{j}^{mod}(x,\bm{\theta})$ {w.r.t.}~$x$ can be approximated by using a finite difference scheme with the already established discretization of size~$h$. A second-order central finite difference scheme is given by
	\begin{align*}
		\frac{\partial^2}{{\partial x}^2} \left( \bm{\phi}_{j}^{mod}(x,\bm{\theta}) \right) \approx \frac{\bm{\phi}_{j}^{mod}(x-h,\bm{\theta}) - 2 \bm{\phi}_{j}^{mod}(x,\bm{\theta}) + \bm{\phi}_{j}^{mod}(x+h,\bm{\theta})}{h^2},
	\end{align*}
	and a second-order forward finite difference scheme is
	\begin{align*}
		\frac{\partial^2}{{\partial x}^2} \left( \bm{\phi}_{j}^{mod}(x,\bm{\theta}) \right) \approx \frac{2\bm{\phi}_{j}^{mod}(x,\bm{\theta})-5\bm{\phi}_{j}^{mod}(x+h,\bm{\theta})+4\bm{\phi}_{j}^{mod}(x+2h,\bm{\theta})-\bm{\phi}_{j}^{mod}(x+3h,\bm{\theta})}{h^2}.
	\end{align*}
	A second-order backward finite difference scheme is defined accordingly. By using these finite difference approximations one obtains the following approximation:
	\begin{align*}
		\frac{\partial^2}{{\partial x}^2} \left( \phi_{j}^{mod}(x,\bm{\theta}) \right) \approx \tunderbrace{\frac{1}{h^2}
			\begin{pmatrix}
				2 & -5 & 4 & -1 & \vphantom{\ddots} & 0 & 0 & 0 & 0 \\
				1 & -2 & 1 & 0 & \ldots\vphantom{\ddots} & 0 & 0 & 0 & 0 \\
				0 & 1 & -2 & 1 & \ddots & 0 & 0 & 0 & 0 \\
				0 & 0 & 1 & -2 & \ddots & 0 & 0 & 0 & 0 \\
				\vdots &\vdots &\ddots & \ddots & \ddots & \ddots & \ddots & \vdots & \vdots \\
				0 & 0 & 0 & 0 &  \ddots & -2 & 1 & 0 & 0 \\
				0 & 0 & 0 & 0 &  \ddots & 1 & -2 & 1 & 0 \\
				0 & 0 & 0 & 0 &  \ldots\vphantom{\ddots} & 0 & 1 & -2 & 1 \\
				0 & 0 & 0 & 0 & \vphantom{\ddots}  & -1 & 4 & -5 & 2
		\end{pmatrix}}{\bm{C}} \begin{pmatrix}
			\bm{\phi}_{j}^{mod}(x_1,\bm{\theta}) \\
			\bm{\phi}_{j}^{mod}(x_2,\bm{\theta}) \\
			\vdots \\
			\bm{\phi}_{j}^{mod}(x_N,\bm{\theta}) \\
		\end{pmatrix}.
	\end{align*}
	Thus, we have
	\begin{align*}
		\frac{\partial^2}{{\partial x}^2} \left( \phi_{j}^{mod}(x,\bm{\theta}) \right) \approx \bm{C} \begin{pmatrix}
			\phi_{j}^{mod}(x_1,\bm{\theta}) \\
			\phi_{j}^{mod}(x_2,\bm{\theta}) \\
			\vdots \\
			\phi_{j}^{mod}(x_N,\bm{\theta}) \\
		\end{pmatrix} = \bm{C} \bm{\phi}_{j}^{mod}(x,\bm{\theta}),
	\end{align*}
	and one also directly obtains
	\begin{align*}
		\frac{\partial^2}{{\partial x}^2} \left( \frac{\partial}{\partial \theta_{k}} \left( \phi_{j}^{mod}(x,\bm{\theta}) \right) \right) \approx \bm{C} \frac{\partial \bm{\phi}_{j}^{mod}(x,\bm{\theta})}{\partial \theta_{k}}.
	\end{align*}
	An evaluation at $x=x_i$ for $i=1,\ldots,N$ yields the $i$-th entry of the resulting vector, and we obtain
	\begin{align*}
		\frac{\partial E_c(\bm{\theta})}{\partial \theta_{k}}  = -2 \sum_{j=1}^{m} \sum_{i=1}^{N} \frac{\kappa_{j}^{exp}(x_i)-\kappa_{j}^{mod}(x_i,\bm{\theta})}
		{\max\left( \left| \kappa_{j}^{exp}(x_i) \right|, \epsilon \right)^2} \cdot \left( \left. \bm{C} \frac{\partial \bm{\phi}_{j}^{mod}(x,\bm{\theta})}{\partial \theta_{k}} \right|_{x=x_i} \right).
	\end{align*}
	
 \bibliographystyle{plainurl}
 \bibliography{literature}

\end{document}